\documentclass[a4,10pt] {article}
\author{Mar\'{\i}a de la Paz Tirado Hern\'andez}
\usepackage{amssymb}
\usepackage{amsfonts}
\usepackage{latexsym}

\usepackage{amsmath}
\usepackage[latin1]{inputenc}

\usepackage[all]{xy}

\usepackage{theorem}
\usepackage{graphicx}
\usepackage{pb-diagram}
\usepackage{verbatim}

\usepackage{anysize}
\usepackage{cancel}
\usepackage{color}
\usepackage[colorlinks=true,linkcolor=blue,urlcolor=black]{hyperref}

\newtheorem{teo}{Theorem}[section]
\newtheorem{prop}[teo]{Proposition}
\newtheorem{cor}[teo]{Corollary}
\newtheorem{lem}[teo]{Lemma}
\newtheorem{Def}[teo]{Definition}
\newtheorem{nota}[teo]{Remark}

\newtheorem{ej}[teo]{Examples}
\newtheorem{hip}[teo]{Hypothesis}
\newtheorem{notacion}[teo]{Notation}

\newcommand{\s}{{\bf s}}
\newcommand{\bt}{{\bf t}}
\newcommand{\I}{\mathcal I}

\DeclareMathOperator{\charr}{\rm char}

\DeclareMathOperator{\Hom}{\rm Hom}

\DeclareMathOperator{\End}{\rm End}

\DeclareMathOperator{\ord}{\rm ord}

\DeclareMathOperator{\Supp}{\rm Supp}

\DeclareMathOperator{\IDer}{\rm IDer}

\DeclareMathOperator{\Der}{\rm Der}

\DeclareMathOperator{\HS}{\rm HS}

\DeclareMathOperator{\Id}{\rm Id}

\usepackage{fancyhdr}
\marginsize{2cm}{2cm}{2.5cm}{2.5cm}

\title{Leaps of modules of integrable derivations in the sense of Hasse-Schmidt}
\author{Mar\'{\i}a de la Paz Tirado Hern\'andez\thanks{Partially supported by MTM2016-75027, P12-FQM-2696 and FEDER.} \thanks{Departamento de \'Algebra e Instituto de Matem\'aticas (IMUS), Universidad de Sevilla, Espa\~na.}}

\begin{document}
\date {}
\maketitle
\begin{abstract}
Let $k$ be a commutative ring of characteristic $p>0$. We prove that
leaps of chain formed by modules of integrable derivations in the
sense of Hasse-Schmidt of a $k$-algebra only occur at powers of $p$.


\end{abstract}

\begin{center} INTRODUCTION \end{center}

Let $k$ be a commutative ring and $A$ a commutative $k$-algebra. A Hasse-Schmidt derivation of $A$ over $k$  of length $m\in \mathbb N$ or $m=\infty$ is a sequence $D=(D_n)_{n\geq 0}^m$ such that: $$
\begin{array}{ccc}
D_0=\Id_A,&\displaystyle D_n(xy)=\sum_{a+b=n} D_a(x)D_b(y)
\end{array}
$$
for all $x,y\in A$. We denote by  $\HS_k(A;m)$ the set of
Hasse-Schmidt derivations of $A$ of length $m$. The component $D_n$ of a Hasse-Schmidt derivation is a differential operator of order $\leq n$ vanishing at 1, in particular $D_1$ is a $k$-derivation.

Hasse-Schmidt derivations of length $m$, also called higher
derivations of order $m$ (see \cite{Ma}), were introduced by H.
Hasse and F.K. Schmidt (\cite{H-S}) and they have been used by
several authors in different contexts (see \cite{Na1},
\cite{HK},\cite{Se}, \cite{Tr} or \cite{Vo}). An important notion
related with Hasse-Schmidt derivations is integrability. Let $m\in
\mathbb N$ or $m=\infty$, then we say that $\delta\in \Der_k(A)$ is
$m$-integrable if there exists $D\in \HS_k(A;m)$ such that
$\delta=D_1$. The set of all $m$-integrable $k$-derivations is an
$A$-submodule of $\Der_k(A)$ for all $m$, which is denoted by
$\IDer_k(A;m)$.

If $k$ has characteristic 0 or $A$ is $0$-smooth over $k$, then any
$k$-derivation is $\infty$-integrable (\cite{Ma}), that means that
$\Der_k(A)=\IDer_k(A;\infty)$. However, if we consider $k$ a ring of
positive characteristic and $A$ any commutative $k$-algebra, then we
do not have the same property, i.e., there could be $n$, a positive
integer, for which $\IDer_k(A;n-1)\neq \IDer_k(A;n)$ (see
\cite{Na2}, \cite{Ti}). In this case, we say that $A$ has a leap at
$n$. Nonetheless, the modules $\IDer_k(A;m)$ have better properties
than $\Der_k(A)$ (see \cite{Mo}) and so their  exploration could
help us to better understand singularities in positive
characteristic.


The aim of this paper is to prove that leaps of modules of
integrable $k$-derivations, where $k$ is a ring of characteristic
$p>0$, only occur at powers of $p$.

This paper is organized as follows:

In section 1 we recall the definition of Hasse-Schmidt derivations
and give some properties that will be useful in later sections.
Moreover, we associate to any Hasse-Schmidt derivation a special
Hasse-Schmidt derivation that we use to prove the main theorem of
this paper.

In section 2 we prove that any $k$-algebra does not have leaps at
certain integers. Namely, if $k$ is any commutative ring and $A$ any
commutative $k$-algebra, we show that $A$ does not have leaps at any
integers invertible in $k$; If the characteristic of $k$ is $p=2$,
then we show that $A$ does not have a leap at 6, and if the
characteristic of $k$ is $p\neq 2$, then we prove that  $A$ does not
have leap at $2p$.

In section 3 we give an integral of the first non-vanishing component of a Hasse-Schmidt derivation and in section 4, we prove our main result, namely that
$A$ only has leaps at powers of $p$.

\begin{center}
Throughout this paper, all rings (and algebras) are assumed to be
commutative.
\end{center}

\section{Hasse-Schmidt derivations}

In this section, we recall the main definitions of the theory of
Hasse-Schmidt derivations and give some results that will be useful
in other sections. From now on, $k$ will be a commutative ring and
$A$ a commutative $k$-algebra. We denote $\overline{\mathbb
N}:=\mathbb N \cup \{\infty\}$ and, for each integer $m\geq 1$, we
will write $A[|\mu|]_m:=A[|\mu|]/\langle \mu^{m+1}\rangle$ and
$A[|\mu|]_\infty:=A[|\mu|]$. General references for the definitions
and results in this section are \cite{Ma} and \cite{Na2}.

\begin{Def}\label{DefHS}
A Hasse-Schmidt derivation of $A$ (over $k$)  of length $m\geq 1$
(resp. of length $\infty$) is a sequence $D:=(D_0,D_1,\ldots, D_m)$
(or resp. $D=(D_0,D_1,\ldots)$) of $k$-linear maps $D_n:A\rightarrow
A$, satisfying the conditions:
$$
\begin{array}{ccc}
D_0=\Id_A,&\displaystyle D_n(xy)=\sum_{i+j=n} D_i(x)D_j(y)
\end{array}
$$
for all $x,y\in A$ and for all $n$. We write $\HS_k(A;m)$ (resp.
$\HS_k(A)$) for the set of Hasse-Schmidt derivations of $A$ (over
$k$) of length $m$ (resp. $\infty$).
\end{Def}

Any Hasse-Schmidt derivation $D\in \HS_k(A;m)$ is determined by the
$k$-algebra homomorphism
$$
\begin{array}{rccc}
\varphi_D:&A&\rightarrow& A[|\mu|]_m \\ &x&\mapsto& \displaystyle
\sum_{n\geq 0}^m D_n(x)\mu^n
\end{array}
$$
satisfying $\varphi_D(x)=x\mod \mu$. $\varphi_D$ can be uniquely
extended to a $k$-algebra automorphism
$\widetilde\varphi_D:A[|\mu|]_m\rightarrow A[|\mu|]_m$ with
$\widetilde\varphi_D(\mu)=\mu$. So, $\HS_k(A;m)$ has a canonical
group structure where the identity is $\Id=(\Id, 0,\ldots,0)$ and
the inverse will be called $D^\ast$. Namely, $D\circ D'=D''\in
\HS_k(A;m)$ with
$$D''_n=\sum_{i+j=n} D_i\circ D_j'$$ for all $n$. Moreover, the $D_n$ component is a $k$-linear differential operator of order $\leq n$ vanishing at 1. In particular, the
 $D_1$ component is a $k$-derivation. So, the map $(\Id, D_1)\in
\HS_k(A;1)\mapsto D_1\in \Der_k(A)$ is an isomorphism of groups.

There are three important operations in this paper:

\begin{enumerate}\label{Tres operaciones con HS}

\item \label{Multporesc} For any $a\in A$ and any $D\in \HS_k(A;m)$, the sequence $a\bullet D=(a^iD_i)\in
\HS_k(A;m)$.

\item For any $1\leq n\leq m$ and any $D\in \HS_k(A;m)$, we define the truncation map by $\tau_{mn}(D)=(\Id,D_1,\ldots,D_n)\in \HS_k(A;n)$.

\item For any $D\in \HS_k(A;m)$ and any $n\geq 1$, we define $D[n]\in \HS_k(A;mn)$ as
    $$
    D[n]_i=\left\{\begin{array}{ll}
    D_{i/n}& \mbox{ if } i=0\mod n\\
    0&\mbox{ otherwise}
    \end{array}\right.
    $$
\end{enumerate}

\begin{Def}\label{Log-Int} Let $D\in \HS_k(A;m)$ where $m\in \overline{\mathbb N}$ and $n\geq m$. Let $I$ be an ideal of $A$.
\begin{itemize}
\item $D$ is $I$-logarithmic if $D_i(I)\subseteq I$ for all $i$. The set of $I$-logarithmic Hasse-Schmidt derivations is denoted by $\HS_k(\log I;m)$, $\HS_k(\log I):=\HS_k(\log I;\infty)$ and
$\Der_k(\log I):=\HS_k(\log I;1)$.
\item  $D$ is $n$-integrable if there exists $E\in \HS_k(A,n)$ such that $\tau_{nm}(E)=D$. Any such $E$
will be called a $n$-integral of $D$. If $D$ is $\infty$-integrable
we say that $D$ is integrable. If $m=1$, we write $\IDer_k(A;n)$ for
the set of $n$-integrable derivations and
$\IDer_k(A):=\IDer_k(A;\infty)$.
\item $D$ is $I$-logarithmically $n$-integrable if there exists $E\in \HS_k(\log I;n)$ such that $E$ is a $n$-integral of
$D$. We put $\IDer_k(\log I;n)$ for the set of $I$-logarithmically
$n$-integrable derivations  when $m=1$ and $\IDer_k(\log
I):=\IDer_k(\log I, \infty)$.
\end{itemize}
\end{Def}

\begin{notacion}\label{Def r-I-log}
Let $D\in \HS_k(A;m)$ be a Hasse-Schmidt derivation where $m\in
\overline{\mathbb N}$, $r< m$ and $I$ an ideal of $A$. We say that
$D$ is $r-I$-logarithmic  if $\tau_{mr}(D)\in \HS_k(\log I;r)$.
 \end{notacion}

\begin{lem}\label{IDer submodulo}
\begin{enumerate}
\item $\HS_k(\log I;m)$ is a subgroup of $\HS_k(A;m)$ for all $m\in \overline{\mathbb N}$.
\item $\IDer_k(A;n)$ is an $A$-submodule of $\Der_k(A)$ for all $n\in \overline{\mathbb N}$ thanks to the group structure of $\HS_k(A;n)$ and operation \ref{Multporesc}.
\end{enumerate}
\end{lem}

\begin{Def}\label{leap}
$A$ has a leap at $s>1$ if the inclusion $\IDer_k(A;s-1)\supsetneq
\IDer_k(A;s)$ is proper.
\end{Def}

\begin{Def}\label{l(D)}
For each Hasse-Schmidt derivation $D\in \HS_k(A;m)$ such that $D\neq
\Id$, we denote $$ \ell(D):=\min\{h\geq 1\mbox{ }|\mbox{ }D_h\neq
0\}
$$
and for $D=\Id$, $\ell(D)=\infty$.
\end{Def}

\begin{lem}[\cite{Na3}, \S 4]\label{l(de la composicion)}
If $D,E\in \HS_k(A;m)$, then $\ell(D\circ E)\geq
\min\{\ell(D),\ell(E)\}$.
\end{lem}

\begin{Def}\label{l(D;e)}
For each $D\in \HS_k(A;m)$ and $1< e\leq m$ integer, if $D_j=0$ for
all $j\neq 0 \mod e$, we denote $\ell(D;e)=\lceil m/e\rceil$ if
$m<\infty$ and $\ell(D;e)=\infty$ if $m=\infty$. Otherwise,
$$
\ell(D;e):=\min\{h\geq 0\mbox{ }|\mbox{ }D_{he+\alpha}\neq 0 \mbox{
for some } \alpha\in \{1,\ldots,e-1\}\}.$$
\end{Def}

\begin{lem}\label{l(D;e) mayor que 0}
\begin{itemize}
\item $\ell(D)\geq e$ if and only if $\ell(D;e)\geq 1$.
\item If $D\in \HS_k(A;m)$ and $1< e\leq m$, then
$\ell(D[e];e)=\lceil m/e\rceil$ if $m<\infty$ and
$\ell(D[e];e)=\infty$ when $m=\infty$.
\item If $\ell(D;e)=i\geq 1$ and $\ell(E;je)\geq i/j$
where $1\leq j\leq i$, then $\ell(D\circ E;e)\geq i$.
\end{itemize}
\end{lem}

\noindent{\bf Proof.} The first two statements are obvious, we will
prove the third one. We denote $D'=D\circ E$. To show that
$\ell(D';e)\geq i$, we have to see that $D'_\alpha=0$ for all
$\alpha <ie$ such that $\alpha\neq 0\mod e$. Let us consider
$\alpha$ with these properties. Since $\ell(E;je)\geq i/j\geq 1$
then $ie\leq \ell(E;je)je$, so we have that $E_\gamma=0$ for all
$\gamma\neq 0 \mod je$ such that $\gamma\leq ie$. Thanks to this,
$$
D_\alpha'=\sum_{\beta+\gamma=\alpha} D_\beta \circ
E_\gamma=\sum_{\gamma=0}^\alpha D_{\alpha-\gamma}\circ
E_\gamma=\sum_{\gamma=0}^{\lfloor\alpha/je \rfloor}
D_{\alpha-je\gamma} \circ E_{je\gamma}
$$
Note that $\alpha-je\gamma\neq 0\mod e$ and $\alpha-je\gamma
<ie-je\gamma\leq ie$. So, $D_{\alpha-je\gamma}=0$
because $\ell(D;e)=i$. Hence, $\ell(D';e)\geq i$.
\begin{flushright}$\square$\end{flushright}

\begin{lem}\label{l(D;e) y derivaciones}
Let $D\in \HS_k(A;m)$ be a Hasse-Schmidt derivation of length $m\in
\overline{\mathbb N}$ and $1< e\leq m$ an integer. Let us assume
that $\ell(D;e)=i\geq 1$, then $D_{ie+\alpha}\in \Der_k(A)$ for all $ie+\alpha\leq m$ where $\alpha=0,\ldots, e-1$.
\end{lem}

\noindent{\bf Proof.} From the definition of Hasse-Schmidt
derivation,
$$
D_{ie+\alpha}(xy)=\sum_{a+b=ie+\alpha} D_a(x)D_b(y)=\sum_{a=0}^{ie}
D_a(x)D_{ie+\alpha-a}(y)+\sum_{a=1}^\alpha
D_{ie+a}(x)D_{\alpha-a}(y).
$$
In the second term, $D_{\alpha-a}=0$ for all $a\neq \alpha$ because
$0<\alpha-a<e$ and $\ell(D;e)\geq 1$. In the first one, since
$\ell(D;e)=i$, if $a\neq 0\mod e$, then $D_a=0$, so we can write the
previous equation as:
$$
D_{ie+\alpha}(xy)=\sum_{a=0}^i D_{ae}(x)D_{ie+\alpha-ae}(y)+
D_{ie+\alpha}(x)y
$$
Note that if $a\neq 0$, then $ie+\alpha-ae<ie$. Moreover
$ie+\alpha-ae\neq 0\mod e$, so $D_{ie+\alpha-ae}=0$. Then,
$$
D_{ie+\alpha}(xy)= xD_{ie+\alpha}(y)+D_{ie+\alpha}(x)y
$$
i.e, $D_{ie+\alpha}$ is a $k$-derivation of $A$ for all
$\alpha=0,\ldots,e-1$.
\begin{flushright}$\square$\end{flushright}

\begin{lem}\label{Integrar una deriv con l(D;m)=n hasta (n+1)m-1}
Let $m>1$ and $n>0$ be two integers and $D\in \HS_k(A;mn)$ a
Hasse-Schmidt derivation such that $\ell(D;m)=n$. Then, $D$ is
$(n+1)m-1$-integrable and there exists an integral of $D$, $D'\in
\HS_k(A;(n+1)m-1)$, such that $\ell(D';m)=n+1$. Moreover, if
$I\subseteq A$ is an ideal and $D\in \HS_k(\log I;mn)$, then $D$ is
$I$-logarithmically $(n+1)m-1$-integrable.
\end{lem}

\noindent{\bf Proof.} Let $\delta_1,\ldots,\delta_{m-1}\in
\Der_k(A)$ be $k$-derivations and let us consider the sequence
$$
D'=(\Id,D'_1,\ldots,D'_{mn},D'_{mn+1},\ldots,D'_{mn+m-1})
=(\Id,D_1,\ldots,D_{mn},\delta_1,\ldots,\delta_{m-1}).
$$  We claim that $D'\in \HS_k(A;(n+1)m-1)$. If this is true, $D'$ is a $(n+1)m-1$-integral of $D$.

To prove this claim we have to show that the following equality must
hold for all $\alpha=1,\ldots, m-1$:
$$
D'_{mn+\alpha}(xy):=\delta_\alpha(xy)=\sum_{\beta=0}^{mn+\alpha}
D'_\beta(x)D'_{mn+\alpha-\beta}(y)
$$
By hypothesis, $D_\beta=0$ for all $\beta\neq 0\mod m$ and
$\beta\leq n$. Since $D'_\beta=D_\beta$ for all $\beta\leq mn$,
$$
\begin{array}{rl}
\displaystyle \sum_{\beta=0}^{mn+\alpha}
D'_\beta(x)D'_{mn+\alpha-\beta}(y)=&
\displaystyle\sum_{\beta=0}^{mn}
D_\beta(x)D'_{mn+\alpha-\beta}(y)+\sum_{\gamma=mn+1}^{mn+\alpha}
D'_\gamma(x)D'_{mn+\alpha-\gamma}(y)\\=&
\displaystyle\sum_{\beta=0}^{n} D_{\beta
m}(x)D'_{(n-\beta)m+\alpha}(y)+\sum_{\gamma=1}^{\alpha}
D'_{mn+\gamma}(x)D'_{\alpha-\gamma}(y)
\end{array}
$$
In the first term, if $\beta>0$, then $0<(n-\beta)m+\alpha<mn$ and
$(n-\beta)m+\alpha\neq 0\mod m$, so
$D'_{(n-\beta)m+\alpha}=D_{(n-\beta)m+\alpha}=0$. In the second one,
if $\gamma\neq \alpha$, then
$D'_{\alpha-\gamma}=D_{\alpha-\gamma}=0$ because
$0<\alpha-\gamma<m$. So,
$$
\sum_{\beta=0}^{mn+\alpha}
D'_\beta(x)D'_{mn+\alpha-\beta}(y)=xD'_{mn+\alpha}(y)+D'_{mn+\alpha}(x)y
=x\delta_\alpha(y)+\delta_\alpha(x)y=\delta_{\alpha}(xy)
$$
Observe that, for each $\alpha=1,\ldots, m-1$, we can choose any
$k$-derivation to be $\delta_\alpha$. In particular, we can put
$\delta_\alpha=0$ for all $\alpha$. In that case, $\ell(D';m)=n+1$.
Thanks to this, it is easy to see that if $D$ is $I$-logarithmic, then $D$ is
$I$-logarithmically $(n+1)m-1$-integrable.
\begin{flushright}$\square$\end{flushright}




\begin{lem}\label{Operacion inversa de desplazar}
Let $m>1$ be an integer and $n\in \overline{\mathbb N}$. If $D\in
\HS_k(A;mn)$ is a Hasse-Schmidt derivation such that $\ell(D;m)=n$
then, there exists $D'\in \HS_k(A;n)$ such that
$D_\alpha'=D_{m\alpha }$ for all $\alpha\leq n$.
\end{lem}

\noindent{\bf Proof.} We have to prove that $D'$ is a Hasse-Schmidt
derivation, so $D'_0=D_0=\Id$ and
$$
\begin{array}{rl}
\displaystyle D_\alpha'(xy)=D_{m\alpha}(xy)=&\displaystyle
\sum_{\beta+\gamma=m\alpha}
D_\beta(x)D_\gamma(y)=\sum_{m\beta+m\gamma=m\alpha}
D_{m\beta}(x)D_{m\gamma}(y)\\ =&\displaystyle
\sum_{\beta+\gamma=\alpha}
D_{m\beta}(x)D_{m\gamma}(y)=\sum_{\beta+\gamma=\alpha}
D'_\beta(x)D'_\gamma(y)
\end{array}
$$
where the third equality holds thanks to $\ell(D;m)=n$.
\begin{flushright}$\square$\end{flushright}

\begin{lem}\label{proced}
Let $D\in \HS_k(A;n)$ be a Hasse-Schmidt derivation of length $n\in
\overline{\mathbb N}$. For each $m>1$, there exists $E\in
\HS_k(A;(n+1)m-1)$ such that $E_m=-D_1$ and $\ell(E;m)=n+1$.
Moreover, if $D$ is $I$-logarithmic for $I\subseteq A$ an ideal,
then $E$ is $I$-logarithmic.
\end{lem}

\noindent{\bf Proof.} We know that $D':=\left( (-1)\bullet
D\right)[m]$ is a Hasse-Schmidt derivation of length $mn$ such that
$D'_m=-D_1$ and $\ell(D';m)=n$. By Lemma \ref{Integrar una deriv con
l(D;m)=n hasta (n+1)m-1}, there exists $E\in  \HS_k(A;(n+1)m-1)$ an
integral of $D'$ with $\ell(E;m)=n+1$. So, this derivation satisfies
the lemma. Moreover, if $D$ is $I$-logarithmic then $D'$ is also
$I$-logarithmic  and, by Lemma \ref{Integrar una deriv con l(D;m)=n
hasta (n+1)m-1}, $E$ is $I$-logarithmic too.
\begin{flushright}$\square$\end{flushright}

\begin{Def}\label{EDb}
For each $D\in \HS_k(A;n)$ and $m>1$, we denote by $E^{D,m}\in
\HS_k(A;(n+1)m-1)$ the Hasse-Schmidt derivation defined in Lemma
\ref{proced}.
\end{Def}

\subsection{Some technical lemmas about composition of Hasse-Schmidt derivations}

In this section, we give some results related with the
composition of Hasse-Schmidt derivations.

\begin{lem}\label{Logaritmico y ultima componente}
Let $I\subseteq A$ be an ideal. If $D\in \HS_k(A;n)$ is
$(n-1)-I$-logarithmic and $E\in \HS_k(\log I;n)$, then $D':=D\circ
E\in \HS_k(A;n)$ is $(n-1)-I$-logarithmic with $D'_n=D_n+H$ where
$H$ is an $I$-logarithmic differential operator, i.e $H(I)\subseteq
I$.
\end{lem}

\noindent{\bf Proof.} From the definition of the composition between
Hasse-Schmidt derivations, we have that
$$
D'_\alpha=\sum_{\beta+\gamma=\alpha} D_\beta \circ
E_\gamma=D_\alpha+\sum_{\beta+\gamma=\alpha,\beta\neq \alpha}
D_\beta\circ E_\gamma
$$
By hypothesis, $D_\beta$ and $E_\gamma$ are $I$-logarithmic for all
$\beta<n$ and all $\gamma\leq n$. Since $\beta+\gamma=\alpha\leq n$
and $\beta\neq \alpha$, the last term is always $I$-logarithmic.
Moreover, if $\alpha<n$, then $D_\alpha$ is $I$-logarithmic, so $D'$
is $(n-1)-I$-logarithmic. On the other hand, if $\alpha=n$,
$D'_n=D_n+H$ where $H=\sum_{\beta+\gamma=\alpha,\beta\neq \alpha}
D_\beta\circ E_\gamma$ which is an $I$-logarithmic differential
operator, so we have the result.
\begin{flushright}$\square$\end{flushright}

\begin{lem}\label{Composicion 1}
Let $e>1$ and $i\geq  1$ be two integers and $n\geq ie$. Let $D,E\in \HS_k(A;n)$ be
two Hasse-Schmidt derivations such that $\ell(D;e)=i\geq 1$ and
$\ell(E)>ie$ and denote $D':=D\circ E\in \HS_k(A;n)$. Then, for
$\alpha \leq n$,
$$
D'_\alpha=\left\{ \begin{array}{ll}
D_\alpha& \alpha\leq ie\\
D_\alpha+E_\alpha& \alpha=ie+1,\ldots, ie+(e-1)
\end{array}\right.
$$
\end{lem}

\noindent{\bf Proof.} If $0<\gamma\leq ie$, then $E_\gamma=0$, so
$$
D'_\alpha=\sum_{\beta+\gamma=\alpha} D_\beta \circ
E_\gamma=D_\alpha+\sum_{\gamma=ie+1}^\alpha D_{\alpha-\gamma}\circ
E_\gamma
$$
Hence, if $\alpha\leq ie$, $D_\alpha'=D_\alpha$. Let us consider
$\alpha=ie+a\leq n$ where $a\in \{1,\ldots, e-1\}$. Then, the previous
equation can be written as
$$
D'_{ie+a}=D_{ie+a}+\sum_{\gamma=1}^a D_{a-\gamma}\circ E_{ie+\gamma}
$$
Note that if $\gamma\neq a$, then $0<a-\gamma<e$ and, since
$\ell(D;e)\geq 1$, $D_{a-\gamma}=0$, i.e,
$D'_{ie+a}=D_{ie+a}+E_{ie+a}$ for all $a$.
\begin{flushright}$\square$\end{flushright}

\begin{lem}\label{composicion 2}
Let $e>1$ and $j>0$ be two integers, $n\geq je$ and $D,E\in \HS_k(A;n)$ two
Hasse-Schmidt derivations such that $\ell(D)=je$ and
$\ell(E;je)=\lceil n/je \rceil$. Let us denote $D':=D\circ E\in
\HS_k(A;n)$. Then, $\ell(D')\geq je$, $\ell(D';e)\geq \ell(D;e)$ and
for each $i\in \mathbb N$ such that $j\leq i\leq \ell(D;e)$, we have
that, for $\alpha\leq n$:
$$
D'_\alpha=\left\{\begin{array}{ll}
D_{je}+E_{je}& \mbox{if } \alpha=je\\
D_\alpha& \mbox{if }\alpha=ie+1,\ldots,ie+e-1
\end{array}
\right.
$$
\end{lem}

\noindent{\bf Proof.} From Lemma \ref{l(de la composicion)},
$l(D')\geq je$. Let us denote $\ell(D;e)=s\geq j$. Then, $(s-1)e<n$,
so $(s-1)/j=(s-1)e/je<\lceil n/je\rceil$. Then, $s-1<\lceil
n/je\rceil j$, i.e, $s/j\leq \ell(E;je)$. Hence,  by Lemma
\ref{l(D;e) mayor que 0}, $\ell(D';e)\geq \ell(D;e)$.

By hypothesis, $E_\gamma=0$ for all $\gamma\neq 0\mod je$ so,
\begin{equation}\label{Equation}
D'_\alpha=\sum_{\beta+\gamma=\alpha} D_\beta \circ
E_\gamma=\sum_{\beta+je\gamma=\alpha} D_\beta\circ E_{je\gamma}
\end{equation}
If $\alpha=je$, then $\gamma$ can only take the values $0$ and $1$,
so $D'_\alpha=D_{je}+E_{je}$. Let us consider $i$ such that  $j\leq
i\leq \ell(D;e)$ and $\alpha=ie+a\leq n$ where $a\in
\{1,\ldots,e-1\}$. Then, in the equation (\ref{Equation}),
$\beta=\alpha-je\gamma=(i-j\gamma)e+a$. Hence, when $\gamma>0$,
$\beta<ie$ and it is not a multiple of $e$, so $D_\beta=0$ and the
only non-zero term is when $\gamma=0$, i.e., $D'_{ie+a}=D_{ie+a}$
for all $a$.
\begin{flushright}$\square$\end{flushright}

The proof of the following lemma is easy by induction:

\begin{lem}\label{Composicion de varias deriv HS}
Let $D^a\in \HS_k(A;n)$ be an orderer family of Hasse-Schmidt
derivations for $a=1,\ldots, t$. We denote $D:=\circ_{a=1}^t
D^a=D^1\circ D^2\circ \cdots \circ D^t\in \HS_k(A;n)$. Then,
$D_\alpha=\sum_{|\beta|=\alpha} D_{\beta_1}^1\circ \cdots \circ
D_{\beta_t}^t$ where $|\beta|=\sum_i \beta_i$.
\end{lem}

\begin{lem}\label{Composicion de derivaciones y l(D;e)}
Let $e,i\geq 1$ be integers and $n\geq ie+e-1$. Let us consider $D^a\in \HS_k(A;n)$ such that
$\ell(D^a;ie+a)\geq 2$ for all $a=1,\ldots,e-2$ and
$\ell(D^{e-1};ie+e-1)\geq 1$. We write $D:=\circ_{a=1}^{e-1}
D^a=D^1\circ D^2\circ \cdots \circ D^{e-1}\in \HS_k(A;n)$. Then,
$\ell(D)\geq ie+1$ and
$$
D_{ie+a}= D^a_{ie+a} \mbox{ where }a=1,\ldots, e-1.
$$
\end{lem}

\noindent{\bf Proof.} Since $\ell(D^a;ie+a)\geq 1$ for all
$a=1,\ldots,e-1$, then $\ell(D^a)\geq ie+a\geq ie+1$ and, by Lemma
\ref{l(de la composicion)}, we can deduce that $\ell(D)\geq ie+1$.
Suppose now that $\alpha=ie+a\leq n$ where $a\in \{1,\ldots, e-1\}$.
From Lemma \ref{Composicion de varias deriv HS}, we have that
$$
D_\alpha= \sum_{|\beta|=\alpha} D_{\beta_1}^1\circ \cdots \circ 
D_{\beta_{e-1}}^{e-1}
$$
Let us consider $\beta=(\beta_1,\ldots,\beta_{e-1})$ such that
$|\beta|=\alpha$. If there is $b\in \{1,\ldots,e-1\}$ such that
$0<\beta_b <ie+b$, then the term associated to $\beta$ is zero so,
we can consider $\beta_b=0$ or $\beta_b\geq ie+b$ for all $b=1,\ldots,e-1$.

Let us suppose that there exist $b,b'\in \{1,\ldots,e-1\}$ such that
$\beta_b,\beta_{b'}>0$, then,
$$
ie+a=\alpha \geq \beta_b+\beta_{b'}\geq ie+b+ie+b'>2ie>\alpha !!!
$$
Hence, there is only one $b\in \{1,\ldots,e-1\}$ such that
$\beta_b\neq 0$. Since $\ell(D^b;ie+b)\geq 2$ for all
$b=1,\ldots,e-2$, we have that $D^b_\gamma=0$ for all
$\gamma=ie+b+1,\ldots,2i(e+b)-1$ (or until $n$ if $n\leq
2i(e+b)-1$). So, in order to the term associated to $\beta$ be not
zero, if  $b\in \{1,\ldots,e-2\}$, $\beta_b=ie+b$ or $\beta_b=0$. On the other
hand, if $b=e-1$ and $\beta_b>ie+b=(i+1)e-1$, then $\alpha=ie+a\leq
(i+1)e-1<\beta_b!!!$ So, $\beta_b=ie+b$ or $\beta_b=0$. Hence, we can conclude
that, if $\beta_b\neq 0$, then $\beta_b=ie+b$ and
$$
\beta_b=ie+b=ie+a=\alpha \Leftrightarrow b=a
$$
Therefore, the only summand which is not zero is the one associated
to $\beta=(0,\ldots,0,ie+a,0,\ldots,0)$ where $ie+a$ is in the
$a$-th position, i.e, $D_{ie+a}=D^a_{ie+a}$ for all
$a=1,\ldots,e-1$.
\begin{flushright}$\square$\end{flushright}

\subsection{Polynomial rings and Hasse-Schmidt derivations}

Let us consider $R=k[x_i|\mbox{ }i \in \mathcal I]$ the polynomial
ring over $k$ in an arbitrary number of variables and $I\subseteq R$
an ideal. In this section, we recall some general result about
integrability of $k$-derivations in polynomial rings.

\begin{teo}{\rm\cite[{\bf Th. 27.1}]{Ma}}\label{IDer=Der}
If the ring $A$ is $0$-smooth over a ring $k$, then a Hasse-Schmidt
derivation of length $m<\infty$ over $k$ can be extended to a
Hasse-Schmidt derivation of length $\infty$.
\end{teo}

\begin{cor}\label{HSenpolise extiende}
Any Hasse-Schmidt derivation of $R$ (over $k$) of length $m\geq 1$
is integrable.
\end{cor}

\noindent{\bf Proof.}
 Since  $R$ is 0-smooth, Theorem \ref{IDer=Der} gives us the corollary.
\vspace{-0.7cm}\begin{flushright}$\square$\end{flushright}




The proof of the following proposition is analogous to that of
Proposition 1.3.4 of \cite{Na2}:

\begin{prop}\label{HS log es sobreyectivo}
Let $R=k[x_i|\mbox{ }i\in \mathcal I]$ be the polynomial ring and
$I\subseteq R$ an ideal. Then, the map $\Pi_n: \HS_k(\log
I;n)\rightarrow \HS_k(R/I;n)$ defined by $\Pi_n(D)=\overline D$
where $\overline D_i(r+I)=D_i(r)+I$ for all $0\leq i\leq n$ is a
surjective group homomorphism.
\end{prop}

The following result generalizes Corollary 2.1.9 of \cite{Na2}.

\begin{cor}\label{Anillo polinomios, logaritmico, map sobrey}
Let $R=k[x_i|\mbox{ }i\in \mathcal I]$ be the polynomial ring and
$I\subseteq R$ an ideal. Then, the map $\Pi: \IDer_k(\log
I;n)\rightarrow \IDer_k(R/I;n)$ defined by $\Pi(\delta)=\overline
\delta$ where $\overline \delta(r+I)=\delta(r)+I$ is a surjective
group homomorphism.
\end{cor}

\noindent{\bf Proof.} Let $\delta\in\IDer_k(R/I;n)$ be a
$n$-integral derivation. From the definition, there exists $D\in
\HS_k(R/I;n)$ an $n$-integral of $\delta$. By Proposition \ref{HS
log es sobreyectivo}, there exists $E\in \HS_k(\log I;n)$ such that
$\Pi_n(E)=D$, in particular $\Pi(E_1)=\delta$ and $E_1\in
\IDer_k(\log I;n)$.
\begin{flushright}$\square$\end{flushright}

\begin{cor}\label{Saltos del cociente}
Let $I$ be an ideal of $R=k[x_i|\mbox{ }i\in \mathcal I]$. Then,
$R/I$ has a leap at $s> 1$ if and only if the inclusion
$\IDer_k(\log I;s-1)\supsetneq \IDer_k(\log I;s)$ is proper.
\end{cor}

\subsection{Multivariate Hasse-Schmidt derivations}

In this section we recall some notions and results of \cite{Na3}.
Throughout this section, $k$ will be a commutative ring and $A$ a
commutative $k$-algebra. Let $q\geq 1$ be an integer and let us call
$\s=\{s_1,\ldots,s_q\}$ a set of $q$ variables.

The monoid $\mathbb N^q$ is endowed with a natural partial ordering.
Namely, for $\alpha,\beta\in \mathbb N^q$, we define
$$
\alpha\leq \beta \Leftrightarrow \exists \gamma \in \mathbb N^q
\mbox{ such that } \beta=\alpha+\gamma \Leftrightarrow \alpha_i\leq
\beta_i \mbox{ } \forall i=1,\ldots,q
$$

The support of a series $a=\sum_\alpha a_\alpha \s^\alpha\in
A[|\s|]$ is $\Supp(a):=\{\alpha\in \mathbb N^q\mbox{ }|\mbox{ }
a_\alpha\neq 0\}$. The order of a non-zero series $a=\sum_\alpha
a_\alpha \s^\alpha\in A[|\s|]$ is
$$
\ord(a):=\min \{|\alpha|\mbox{ }|\mbox{ } \alpha\in \Supp(a)\}
$$
and if $a=0$ we define $\ord(a):=\infty$.

\begin{Def}\label{co-ideal}
We say that a subset $\Delta\in \mathbb N^q$ is a co-ideal of
$\mathbb N^q$  if whenever $\alpha\in \Delta$ and $\alpha'\leq
\alpha$, then $\alpha'\in \Delta$.
\end{Def}

For example, for $\beta\in \mathbb N^q$,  $\mathfrak
n_\beta:=\{\alpha\in \mathbb N^q\mbox{ }| \mbox{ } \alpha\leq
\beta\}$ is a co-ideal of $\mathbb N^q$.

\begin{Def}\label{DeltaA}
For each co-ideal $\Delta\subset \mathbb N^q$, we  denote by
$\Delta_A$ the ideal of $A[|\s|]$ whose elements are the series
$\sum_{\alpha\in \mathbb N^q} a_\alpha \s^\alpha$ such that
$a_\alpha=0$ if $\alpha\in \Delta$. i.e, $\Delta_A=\{a\in A[|\s|]
\mbox{ }|\mbox{ } \Supp(a)\subseteq \Delta^c\}$.
\end{Def}

Let us denote $A[|\s|]_\Delta:=A[|s|]/\Delta_A$. Note that if $q=1$
and $\Delta=\{i\mbox{ }|\mbox{ }i\leq m\}$, then
$A[|\s|]_\Delta=A[|s|]_m$ defined before. From now on, $\Delta$ will
be a non-empty co-ideal.

\begin{Def}\label{Def multiderivacion}
A $(q,\Delta)$-variate Hasse-Schmidt derivation of $A$ over $k$ is a
family $D=(D_\alpha)_{\alpha\in\Delta}$ of $k$-linear maps
$D_\alpha:A\rightarrow A$, satisfying the conditions:
$$
D_0=\Id_A, \mbox{   } D_\alpha(xy)=\sum_{\beta+\gamma=\alpha}
D_\beta(x)D_\gamma(y)
$$
for all $x,y\in A$ and for all $\alpha\in \Delta$. We denote by
$\HS_k^{q}(A;\Delta)$ the set of all $(q,\Delta)$-variate
Hasse-Schmidt derivations of $A$ over $k$ and $\HS_k^{q}(A)$ for
$\Delta=\mathbb N^q$. For $q=1$ and $\Delta=\{i|\mbox{ } i\leq m\}$, a
$(1,\Delta)$-variate Hasse-Schmidt derivation is a Hasse-Schmidt
derivation of length $m$ in the  usual way.
\end{Def}

\begin{nota}[\cite{Na3}, \S 4,7]\label{HS en End}
Any $(q,\Delta)$-variate Hasse-Schmidt derivation $D$ of $A$ over
$k$ can be understood as a power series
$$
\sum_{\alpha\in \Delta} D_\alpha \s^\alpha\subseteq
\End_k(A)[|\s|]_\Delta
$$
and so we can consider $\HS_k^q(A;\Delta)\subseteq
\End_k(A)[|\s|]_\Delta$.
\end{nota}

\begin{cor}[\cite{Na3}, Corollary 1]\label{Grupo}
Let $k$ be a ring, $A$ a $k$-algebra, $q\geq
1$ an integer and $\Delta\subseteq \mathbb N^q$ a non-empty
co-ideal. Then, $\HS_k^q(A;\Delta)$ is a group.
\end{cor}

Namely, the group operation in $\HS_k^q(A;\Delta)$ is explicitly
given by
$$
(D,E)\in \HS_k^{q}(A;\Delta)\times \HS_k^{q}(A;\Delta)\mapsto D\circ
E\in \HS_k^q(A;\Delta)
$$
with
$$
(D\circ E)_\alpha=\sum_{\beta+\gamma=\alpha} D_\beta \circ E_\gamma
$$

\begin{notacion}\label{Hom0k-alg}
Let us denote
$$
\Hom_{k-alg}^\circ (A;A[|\s|]_\Delta):=\{f\in
\Hom_{k-alg}(A,A[|\s|]_\Delta)\mbox{ }|\mbox{ } f(x)\equiv x\mod
(\mathfrak n_0)_A \mbox{ } \forall x\in A\}.
$$
\end{notacion}

\begin{lem}[\cite{Na3}, \S 4]\label{Isomor grupo en multideriv}
Let $k$ be a ring, $A$ a $k$-algebra, $q\geq
1$ an integer, $\s=\{s_1,\ldots,s_q\}$ a set of $q$ variables and
$\Delta$ a non-empty co-ideal. Then, the map
$$
D\in \HS_k^q(A;\Delta)\mapsto \left[x\in A\mapsto \sum_{\alpha\in
\Delta} D_\alpha(x)\s^\alpha\right]\in \Hom_{k-akg}^\circ
\left(A,A[|\s|]_\Delta\right)
$$
 is a group isomorphism.
\end{lem}

\subsubsection{Substitutions}

Let $k$ be a commutative ring, $A$ a commutative $k$-algebra,
$\s=\{s_1,\ldots,s_q\}$, $\bt=\{t_1,\ldots,t_m\}$ two sets of
variables where $q,m\geq 1$ and $\Delta \subseteq \mathbb N^q$ and
$\nabla \subseteq \mathbb N^m$ non-empty co-ideals.

\begin{Def}\label{Def sustitucion}
An $A$-algebra map $\phi:A[|\s|]_\Delta\rightarrow A[|\bt|]_\nabla$
will be called a substitution map if $\ord(\phi(s_i))\geq 1$ for all
$i=1,\ldots,q$.
\end{Def}

\begin{prop}[\cite{Na3}, Prop. 10]\label{Sustituciones y Derivaciones de
HS} For any substitution map $\phi:A[|\s|]_\Delta \rightarrow
A[|{\bf t}|]_\nabla$, we have that if $f\in \Hom_{k-alg}^\circ
(A,A[|\s|]_\Delta)$, then $\phi \circ f\in \Hom_{k-alg}^\circ
(A,A[|{\bf t}|]_\nabla)$.
\end{prop}

\begin{notacion} Let $\phi:A[|\s|]_\Delta \rightarrow
A[|{\bf t}|]_\nabla$ be a substitution map and
$\varphi_D:A\rightarrow A[|\s|]_\Delta\in
\Hom^0_{k-alg}(A,A[|\s|]_\Delta)$ the $k$-algebra homomorphism
associated to $D\in \HS_k^q(A;\Delta)$. We denote by $\phi\bullet
D\in \HS_k^m(A;\nabla)$ the $(m,\nabla)$-Hasse-Schmidt derivation
associated to $\phi\circ \varphi_D$.
\end{notacion}

Let $\phi:A[|\s|]_\Delta \rightarrow A[|{\bf t}|]_\nabla$ be a
substitution map and $D=\sum D_\alpha \s^\alpha\in
\HS_k^q(A,\Delta)$, then
$$
\phi \bullet D=\phi\left(\sum D_\alpha \s^\alpha\right)=\sum\phi(\s)^\alpha
D_\alpha
$$

\begin{nota}\label{Estable por sustituciones}
Thanks to the previous expression, it is easy to see that, if
$\phi:A[|\s|]_\Delta\rightarrow A[|\bf t|]_\nabla$ is a substitution
map and $D\in \HS_k^q(\log I;\Delta)$ for any $I\subseteq A$ an
ideal, i.e, $D_\alpha(I)\subset I$ for all $\alpha\in \Delta$, then
$\phi\bullet D\in \HS_k^m(\log I;\nabla)$.
\end{nota}

\begin{ej}\label{Ejemplos de sustituciones} The operations defined in \ref{Tres operaciones con
HS} are examples of substitution maps. Namely, let $D\in \HS_k(A;m)$ a
Hasse-Schmidt derivation of length $m\in \overline{\mathbb N}$.
\begin{enumerate}
\item For any $a\in A$, $a\bullet D=\phi \bullet D$ where
$\phi:\mu\in A[|\mu|]_m\mapsto a\mu\in A[|\mu|]_m$.
\item Let $1\leq n\leq m$ be an integer. If $\phi:\mu\in
A[|\mu|]_m\mapsto \mu\in A[|\mu|]_n$ then, $\tau_{mn}(D)=\phi\bullet
D$.
\item For any $n\geq 1$, $D[n]=\phi\bullet D$ where
$\phi:\mu\in A[|\mu|]_m\mapsto \mu^n\in A[|\mu|]_{mn}$.
\end{enumerate}
\end{ej}

\begin{notacion}\label{Def E}
Let $D\in \HS_k(A)$ be a Hasse-Schmidt derivation. We denote
$B^D=\phi \bullet D\in \HS_k^2(A)$ where $\phi:\mu\in
A[|\mu|]\mapsto \mu_1+\mu_2\in A[|\mu_1,\mu_2|]$.
\end{notacion}

\begin{lem}\label{Componentes de E}
Let $D\in \HS_k(A)$ be a Hasse-Schmidt derivation. Then,
$B^D_{(i,j)}=\binom{i+j}{i} D_{i+j}$ for all $(i,j)\in \mathbb N^2$.
\end{lem}

\noindent{\bf Proof.} We can write $D=\sum_{\alpha\geq 0}
D_\alpha\mu^\alpha\subseteq \End_k(A)[|\mu|]$. Then,
$$
B^D=\phi\bullet\left(\sum_{\alpha\geq 0}
D_\alpha\mu^\alpha\right)=\sum_{\alpha\geq 0}
D_\alpha(\mu_1+\mu_2)^\alpha=\sum_{\alpha\geq
0}D_{\alpha}\sum_{i+j=\alpha}\binom{\alpha}{j}\mu_1^i\mu_2^j=
\sum_{i+j\geq 0}\binom{i+j}{j}D_{i+j}\mu_1^i\mu_2^j
$$
So,
$$
B^D_{(i,j)}=\binom{i+j}{j}D_{i+j}
$$
\begin{flushright}$\square$\end{flushright}

\begin{lem}\label{I-logaritmico para derivacion E}
Let $I$ be an ideal of $A$ and let us consider $D\in \HS_k(A)$ a
$(n-1)-I$-logarithmic Hasse-Schmidt derivation. If $i+j<n$, then
$B^D_{(i,j)}(I)\subseteq I$.
\end{lem}

\noindent{\bf Proof.} If $i+j<n$, then $D_{i+j}(I)\subseteq I$, so
$B^D_{(i,j)}(I)=\binom{i+j}{i}D_{i+j}(I)\subseteq I$.
\vspace{-0.7cm}\begin{flushright}$\square$\end{flushright}

\subsubsection{External product}

\begin{Def}\label{Producto exterior}
Let $R$ be a ring, $q,m\geq 1$, $\s=\{s_1,\ldots,s_q\}$, ${\bf
t}=\{t_1,\ldots,t_m\}$ disjoint sets of variables and $\Delta\subset
\mathbb N^q$ and $\nabla\subset \mathbb N^m$ non-empty co-ideals.
For each $r\in R[|\s|]_\Delta$, $r'\in R[|{\bf t}|]_\nabla$, the
external product $r\boxtimes r'\in R[|\s \sqcup {\bf
t}|]_{\Delta\times \nabla}$ is defined as
$$
r\boxtimes r':=\sum_{(\alpha,\beta)\in \Delta\times \nabla} r_\alpha
r'_{\beta} \s^\alpha {\bf t}^\beta
$$
\end{Def}

\begin{prop}[\cite{Na3}, Prop. 6]\label{Proposicion 6 Na3}
Let $D\in \HS_k^q(A;\Delta)$, $E\in \HS_k^m(A;\nabla)$ be
Hasse-Schmidt derivations. Then its external product $D\boxtimes E$
is a $(\s \sqcup {\bt},\nabla \times \Delta)$-variate Hasse-Schmidt
derivation.
\end{prop}

\begin{nota}\label{Componentes del producto exterior}
With the above notation, $(D\boxtimes E)_{(i,j)}=D_iE_j$ for all
$(i,j)\in \mathbb N^2$.
\end{nota}

\begin{notacion}\label{Def F}
Let $D\in \HS_k(A)$ be a Hasse-Schmidt derivation. We denote
$F^D=D\boxtimes D\in \HS_k^2(A)$ and $\left(F^D\right)^\ast\in
\HS_k^2(A)$ its inverse.
\end{notacion}

It is easy to proof the next lemma:

\begin{lem}\label{Inversa de F}
$\left(F^D\right)^\ast_{(i,j)}=D^\ast_jD^\ast_i$ where $D^\ast\in
\HS_k(A)$ is the inverse of $D$.
\end{lem}

\begin{lem}\label{I-logaritmico para derivacion F}
Let $I$ be an ideal of $A$ and let us consider $D\in \HS_k(A)$ a
$(n-1)-I$-logarithmic Hasse-Schmidt derivation. If $i,j<n$, then
$\left(F^D\right)_{(i,j)}^\ast(I)\subseteq I$.
\end{lem}

\noindent{\bf Proof.} Since $D$ is $(n-1)-I$-logarithmic, $D^\ast$
is $(n-1)-I$-logarithmic too. So,
$\left(F^D\right)^\ast_{(i,j)}(I)=D_j^\ast D_i^\ast(I)\subseteq I$.
\begin{flushright}$\square$\end{flushright}

\subsection{A special Hasse-Schmidt derivation}

In this section, we define a Hasse-Schmidt derivation that we
will use in later sections and give some properties about it. Throughout
this section, $k$ will be a commutative ring, $A$ a commutative
$k$-algebra, $I\subseteq A$ an ideal and $D\in \HS_k(A)$ a
Hasse-Schmidt derivation.

\begin{notacion}\label{Definicion de G}
For each $D\in \HS_k(A)$, we define $G^D:=B^D\circ (F^D)^\ast\in
\HS_k^2(A)$  (see Notations \ref{Def E} and \ref{Def F}). Namely,
$G^D_{(i,j)}=\sum_{\alpha+\beta=(i,j)} B^D_\alpha\circ
(F^D)^\ast_\beta$.
\end{notacion}

From now on, we will omit the superscript and we will write
$G:=G^D$, $B:=B^D$ and $F:=F^D$.

\begin{lem}\label{G(m,0)=G(0,m)=0}
For each $m>0$, we have that $G_{(m,0)}=G_{(0,m)}=0$ and
$G_{(1,m)},G_{(m,1)}\in \Der_k(A)$.
\end{lem}

\noindent{\bf Proof.} First, we calculate $G_{(m,0)}$:
$$
G_{(m,0)}=\sum_{\alpha+\beta=(m,0)} B_\alpha
F^\ast_\beta=\sum_{\alpha_1+\beta_1=m} B_{(\alpha_1,0)}
F^\ast_{(\beta_1,0)}=\sum_{\alpha_1+\beta_1=m} D_{\alpha_1}
D^\ast_{\beta_1}=0
$$
The calculation of $G_{(0,m)}$ is analogous. Now, by definition of
multivariate Hasse-Schmidt derivation: $$
\begin{array}{rl}
G_{(1,m)}(xy)&\displaystyle =\sum_{\substack{\alpha_1+\beta_1=1\\
\alpha_2+\beta_2=m}}
G_{(\alpha_1,\alpha_2)}(x)G_{(\beta_1,\beta_2)}(y)=
\sum_{\alpha_2+\beta_2=m} G_{(0,\alpha_2)}(x)
G_{(1,\beta_2)}(y)+\sum_{\alpha_2+\beta_2=m}
G_{(1,\alpha_2)}(x)G_{(0,\beta_2)}(y) \\&=
xG_{(1,m)}(y)+G_{(1,m)}(x)y
\end{array}
$$
It is analogous for $G_{(m,1)}$.
\begin{flushright}$\square$\end{flushright}

\begin{lem}\label{I-logaritmico para derivacion G} Let us suppose that $D\in \HS_k(A)$ is $(n-1)-I$-logarithmic. We have the following properties:
\begin{enumerate}
\item If $0\leq i+j<n$, then $G_{(i,j)}(I)\subseteq I$
\item If $i$ and $j$ are not zero and $i+j=n>0$, then $G_{(i,j)}=\binom{n}{i}D_n+H$ where $H$ is an $I$-logarithmic differential operator.
\end{enumerate}
\end{lem}

\noindent{\bf Proof.}
\begin{enumerate}
\item If $i+j=0$, then $G_{(i,j)}=\Id$ and, if $i=0$ or $j=0$ then, $G_{(i,j)}=0$ so the result is obvious and we can suppose that $i,j>0$. We have that
$$
    G_{(i,j)}=\sum_{\substack{\alpha_1+\beta_1=i\\ \alpha_2+\beta_2=j}}
B_{(\alpha_1,\alpha_2)} \circ F^\ast_{(\beta_1,\beta_2)}
$$
Since $i$ and $j$ are not zero, $1\leq i,j<n-1$ so,
$\beta_1,\beta_2<n-1$. Moreover,
$\alpha_1+\beta_1+\alpha_2+\beta_2=i+j<n$, so $\alpha_1+\alpha_2<n$.
By Lemmas \ref{I-logaritmico para derivacion E} and
\ref{I-logaritmico para derivacion F}, the sum is $I$-logarithmic.

\item By definition,
$$
\begin{array}{rl}
  \displaystyle  G_{(i,j)}=&\displaystyle \sum_{\substack{ \alpha_1+\beta_1=i\\ \alpha_2+\beta_2=j}}
B_{(\alpha_1,\alpha_2)} \circ F^\ast_{(\beta_1,\beta_2)}=
B_{(i,j)}+\sum_{\substack{\alpha_1+\beta_1=i\\ \alpha_2+\beta_2=j\\
\alpha\neq (i,j)}} B_{(\alpha_1,\alpha_2)} \circ
F^\ast_{(\beta_1,\beta_s)}
\\& \displaystyle
=\binom{n}{i}D_n+\sum_{\substack{\alpha_1+\beta_1=i\\
\alpha_2+\beta_2=j\\ \alpha\neq (i,j)}} B_{(\alpha_1,\alpha_2)}
\circ F^\ast_{(\beta_1,\beta_2)}
\end{array}
$$
If $\alpha\neq (i,j)$, then $\alpha_1<i$ or $\alpha_2<j$ so,
$\alpha_1+\alpha_2<i+j=n$ and, by Lemma \ref{I-logaritmico para
derivacion E}, $B_\alpha(I)\subseteq I$. On the other hand,
$\beta_1,\beta_2<n$ because $i,j<n$. Hence,
$F^\ast_\beta(I)\subseteq I$ (Lemma \ref{I-logaritmico para
derivacion F}). So, the sum is an $I$-logarithmic  differential operator.
\end{enumerate}
\begin{flushright}$\square$\end{flushright}

From now on, $k$ will be a commutative ring of characteristic $p>0$,
$A$ and $I$ as before and $n=e_sp^s+\cdots+e_tp^t$ a positive
integer expressed in base $p$ expansion where $s\geq t\geq 1$ and
$0\leq e_i<p$ with $e_s,e_t\neq 0$ (note that $t$ and $s$ could be
equal). It is easy to proof the next lemma:

\begin{lem}\label{menor con binomio no cero}
Let $p,n$ be as before. Then,
$$
p^t=\min\left\{m\in \mathbb N_+ \mbox{ }|\mbox{ } \binom{n}{m}\neq 0
\mod p\right\}.
$$
\end{lem}

Thanks to this lemma, we can prove the next result:

\begin{lem}\label{G(i,j) con i+j=n i<pt}
Let $p$ be a prime and $n=e_sp^s+\cdots+e_tp^t$ a positive integer
expressed in base $p$ expansion where $e_s,e_t\neq 0$ and $s\geq
t\geq 1$. Let us consider $i,j\geq 0$ such that $i+j=n$ and $i<p^t$.
If $D\in \HS_k(A)$ is $(n-1)-I$-logarithmic then,
$G_{(i,j)}(I)\subseteq I$.
\end{lem}

\noindent{\bf Proof.}  By Lemma \ref{G(m,0)=G(0,m)=0}, if $i=0$ or
$j=0$, then $G_{(i,j)}=0$ so, it is $I$-logarithmic. If $i,j\geq 1$,
by Lemma \ref{I-logaritmico para derivacion G},
$G_{(i,j)}=\binom{n}{i}D_n+H$ where $H(I)\subseteq I$. By Lemma
\ref{menor con binomio no cero}, $\binom{n}{i}=0$ and we have the
result.
\begin{flushright}$\square$\end{flushright}

Let us consider the following substitution map:
$$
\begin{array}{rccl}
\varphi^{r}:&R[|\mu_1,\mu_2|]&\rightarrow& R [|\mu|]\\
&\mu_1&\mapsto & \mu^{r+1}\\
&\mu_2&\mapsto & \mu^r
\end{array}
$$

\begin{notacion}\label{Definicion de Gpt}
Let $p$ be a prime and $n=e_sp^s+\cdots+e_tp^t$ a positive integer
expressed in base $p$ expansion where $s\geq t\geq 1$ and $0\leq
e_i<p$ with $e_s,e_t\neq 0$. Let $D\in \HS_k(A)$ be a Hasse-Schmidt
derivation and let us consider $G^D\in \HS_k^2(A)$ defined in
\ref{Definicion de G}. We define
$G^{D,p^t}=\tau_{\infty,(n+1)p^t}\left(\varphi^{p^t}\bullet
G^D\right)\in \HS_k(A;(n+1)p^t)$.
\end{notacion}

\begin{lem}\label{I-logaritmico para derivacion Gpt}
Let $p,n$ be two positive integers as before. Then,
$\ell\left(G^{D,p^t}\right)\geq 2p^t+1$. Moreover, if $D\in
\HS_k(A)$ is $(n-1)-I$-logarithmic then, $G^{D,p^t}$ is
$((n+1)p^t-1)-I$-logarithmic and
$G^{p^t}_{(n+1)p^t}=\binom{n}{p^t}D_{n}+H$ where $H$ is an
$I$-logarithmic differential operator.
\end{lem}

\noindent{\bf Proof.} Note that
$$
\varphi^{p^t}\bullet G^D=\varphi^{p^t}\left(\sum_{(i,j)}
G^D_{(i,j)}\mu_1^i\mu_2^j\right)=\sum_{(i,j)}
G^D_{(i,j)}\mu^{(p^t+1)i+p^tj}=\sum_{\alpha\geq 0}\left(
\sum_{(i,j):(p^t+1)i+p^tj=\alpha} G^D_{(i,j)}\right)\mu^\alpha
$$
Since $G^D_{(i,j)}=0$ if $i$ or $j$ is zero (Lemma
\ref{G(m,0)=G(0,m)=0}), we have that $G^{D,p^t}_0=\Id$ and for all
$\alpha\geq 1$,
$$
G^{D,p^t}_\alpha=\sum_{\substack{(i,j):(p^t+1)i+p^tj=\alpha\\
i,j\neq 0}} G^D_{(i,j)}
$$
If $\alpha<2p^t+1$ then there is  not $(i,j)$ with $i,j\neq 0$ such
that $(p^t+1)i+p^tj=\alpha$, so $G^{D,p^t}_\alpha=0$. Hence,
$\ell\left(G^{D,p^t}\right)\geq 2p^t+1$. Now, we will suppose that
$D$ is $(n-1)-I$-logarithmic and will prove the rest of the lemma.

Let us consider a pair $(i,j)$ with $i,j\neq 0$ and $i+j=n+l$ where
$l\geq 0$. Then,
$$
(p^t+1)i+p^tj=p^t(i+j)+i=p^t(n+l)+i.
$$

If $l>0$, then $p^t(n+l)+i>p^t(n+l)\geq p^t(n+1)$. So, $G^D_{(i,j)}$
does not appear in any component of $G^{D,p^t}$.

If $l=0$, then $p^tn+i\leq (n+1)p^t$ if and only if $i\leq p^t$. So,
$G^D_{(i,j)}$ appears in some component of $G^{D,p^t}$ if $i\leq
p^t$. By Lemma \ref{G(i,j) con i+j=n i<pt}, $G^D_{(i,j)}(I)\subseteq
I$ if $i<p^t$. On the other hand, if $i=p^t$, then $j=n-p^t$ and
$(p^t+1)p^t+p^t(n-p^t)=(n+1)p^t$. Hence, $G^D_{(p^t,n-p^t)}$ is a
term of $G^{D,p^t}_{(n+1)p^t}$ and it is the only component that is
not $I$-logarithmic. So, $G^{D,p^t}$ is $((n+1)p^t-1)-I$-logarithmic
and
$$
G^{D,p^t}_{(n+1)p^t}= G^D_{(p^t,n-p^t)}+\mbox{\it Some
$I$-logarithmic differential operator}=\binom{n}{p^t}D_n+\mbox{\it
Some $I$-logarithmic diff. op.}
$$
where the last equality holds because of Lemma \ref{I-logaritmico para
derivacion G}.
\begin{flushright}$\square$\end{flushright}

\section{Some partial integrability results}

In this section, $k$ will be a commutative ring of characteristic
$p>0$ and $A$ a commutative $k$-algebra. We will give some results
about leaps of modules of integrable $k$-derivations of $A$.
Namely, we prove that $A$ does not have leaps at the integers that
are not a multiple of $p$ and on the first multiple of $p$ which is
not a power of $p$.

\begin{lem}\label{M invertible en k}
If $m$ is invertible in $k$, any Hasse-Schmidt derivation of length
$m-1$ is $m$-integrable.
\end{lem}

\noindent{\bf Proof.} Since $A$ is a $k$-algebra, we can write
$A:=R/I$ where $R$ is a polynomial ring (in an arbitrary number of
variables) and $I\subseteq R$ an ideal. Let $D\in \HS_k(A;m-1)$ be a
Hasse-Schmidt derivation of $A$ of length $m-1$. Then, there exists
$\widetilde D\in \HS_k(\log I;m-1)$ such that $\Pi_{m-1}(\widetilde
D)=D$. Thanks to Corollary \ref{HSenpolise extiende}, we can
integrate $\widetilde D$, so we have $E\in \HS_k(R;m)$ such that
$\tau_{m,m-1}(E)=\widetilde D$. From Definition 1.2.11 and
Proposition 3.1.2 of \cite{Na4}, $\varepsilon_m(E)=mE_m+H\in
\Der_k(R)$ where $H$ is an $I$-logarithmic differential operator.
Then,
$$
E':=E\circ ((-1/m)\bullet
(\Id,\varepsilon_m(E)))[m]=(\Id,E_1,\ldots,E_{m-1},-(1/m)H)\in
\HS_k(\log I;m).
$$
So, $\Pi_m(E')\in \HS_k(A;m)$ is a $m$-integral of $D$ (Proposition
\ref{HS log es sobreyectivo}).
\begin{flushright}$\square$\end{flushright}

\begin{cor}\label{No tiene saltos en los no multiplos de p}
If $k$ has characteristic $p>0$ and $m\neq 0\mod p$. Then,
$\IDer_k(A;m-1)=\IDer_k(A;m)$, i.e., $A$ does not have a leap at $m$.
\end{cor}

\begin{prop}\label{Caso p=2}
Let $k$ be a ring of characteristic $p=2$ and $A$ a
 $k$-algebra. Then, $\IDer_k(A;5)=\IDer_k(A,6)$.
\end{prop}

\noindent{\bf Proof.}  As in the previous proof, we can write
$A:=R/I$ where $R$ is a polynomial ring and $I\subseteq R$ an ideal.
By Corollary \ref{Saltos del cociente}, $\IDer_k(A;5)=\IDer_k(A;6)$
if and only if $\IDer_k(\log I;5)=\IDer_k(\log I,6)$. The inclusion
$\IDer_k(\log I;6)\subseteq \IDer_k(\log I;5)$ is always true, so
let $\delta\in \IDer_k(\log I;5)$ be an $I$-logarithmically
$5$-integrable $k$-derivation and we consider $D\in \HS_k(\log I;5)$
an integral of $\delta$. By Corollary \ref{HSenpolise extiende}, we
can integrate $D$ until $\infty$. So, we have $D=(\Id,
D_1,\ldots,D_5,D_6,\ldots)\in \HS_k(R)$ which is $5-I$-logarithmic.
Then, let us consider $G:=G^D\in \HS_k^2(R)$ defined in
\ref{Definicion de G}. By Lemma \ref{I-logaritmico para derivacion
G}, $G_{(i,j)}(I)\subseteq I$ for all $i+j\leq 5$. Moreover,
$G_{(2,4)}=\binom{6}{2}D_6 +H=D_6+H$ where $H(I)\subseteq I$.

On the other hand, by definition of multivariate Hasse-Schmidt
derivation and Lemma \ref{G(m,0)=G(0,m)=0}:
$$
\begin{array}{rl}
G_{(2,4)}(xy)&\displaystyle =\sum_{\alpha+\beta=(2,4)}
G_\alpha(x)G_\beta(y)= \sum_{\substack{\alpha_1+\beta_1=2\\
\alpha_2+\beta_2=4}}
G_{(\alpha_1,\alpha_2)}(x)G_{(\beta_1,\beta_2)}(y) \\ &\displaystyle
=G_{(2,4)}(x)y +xG_{(2,4)}(y)+
\sum_{\alpha_2+\beta_2=4} G_{(1,\alpha_2)}(x)G_{(1,\beta_2)}(y)\\[0.6cm]
\displaystyle &=G_{(2,4)}(x)y
+xG_{(2,4)}(y)+G_{(1,1)}(x)G_{(1,3)}(y)+
G_{(1,3)}(x)G_{(1,1)}(y)+G_{(1,2)}(x)G_{(1,2)}(y)
\end{array}$$
Since $G_{(1,j)}\in \Der_k(A)$ by Lemma \ref{G(m,0)=G(0,m)=0},
$$
D'=(\Id, G_{(1,2)},G_{(2,4)}-G_{(1,1)}G_{(1,3)})\in \HS_k(R;2)
$$
and $D'$ is $1-I$-logarithmic. Moreover,
$\left(G_{(1,1)}G_{(1,3)}\right)(I)\subseteq I$, so $D'_2=D_6+H'$
where $H'(I)\subseteq I$. Then,
$$
D''=\tau_{\infty,6}(D) \circ D'[3]=(\Id,D_1,\ldots,
D_6+D_3G_{(1,2)}+D_6+H')=(\Id,D_1,\ldots,D_3 G_{(1,2)}+H')\in
\HS_k(\log I;6)
$$
Hence, $\IDer_k(\log I;5)=\IDer_k(\log I;6)$ and we have the result.
\begin{flushright}$\square$\end{flushright}

Now, we prove that $\IDer_k(A;2p-1)=\IDer_k(A;2p)$ when $p\neq 2$.
We will start with some previous results.

\begin{Def}\label{spn}
Let $p$ be a prime and $n=e_sp^s+\cdots +e_0$ a positive integer
expressed in base $p$ expansion where $e_s\neq0$. We define
$s_p(n):=\sum_{i=0}^s e_i$.
\end{Def}

\begin{nota}\label{Notas acerca de spn}
If $1\leq n\leq p-1$, then $s_p(n)=n$. If $n\geq p$, then
$s_p(n)<n$.
\end{nota}

\begin{Def}\label{spn recursivo}
For each $j\geq 0$, we define
$s_p^j(n):=\underbrace{s_p(s_p(\cdots(s_p}_{j \mbox{
times}}(n))\cdots)$.
\end{Def}

\begin{lem}\label{Existencia de Tpn}
There exists $j\geq 0$ such that $s_p^j(n)=s_p^{j+1}(n)$. Moreover,
if $s_p^j(n)=s_p^{j+1}(n)$ then, $s_p^j(n)=s_p^J(n)$ for all $J\geq
j$.
\end{lem}

\noindent{\bf Proof.} If $n\leq p-1$,  $n=s_p(n)$. Hence, the lemma holds
for $j=0$.
If $n\geq p$, then $s_p(n)<n$.
So, if $s_p(n)\leq p-1$, then $s_p^2(n)=s_p(n)$ and the lemma holds
for $j=1$. Otherwise, $s_p^2(n)<s_p(n)<n$. By performing this process
recursively, we obtain that $s_p^j(n)\leq p-1$ for some $j$. So,
$s_p^j(n)=s_p^{j+1}(n)$  and the lemma holds for this $j$. Moreover, if $s_p^j(n)=s_p^{j+1}(n)$, then
$s_p^j(n)\leq p-1$, so $s_p^{J}(n)=s_p^j(n)$ for all $J\geq j$.
\begin{flushright}$\square$\end{flushright}

\begin{Def}\label{Tpn}
Let $p$ be a prime and $n$ a positive integer. Let us consider
$j=\min\{l\geq 0\mbox{ } |\mbox{ } s_p^l(n)=s_p^{l+1}(n)\}$. We
define $T_p(n):=s_p^j(n)$.
\end{Def}

\begin{lem}\label{Igualdad de Tpn}
$T_p(n)=s_p^J(n)$ for all $J\in\{l\geq 0\mbox{ } |\mbox{ }
s_p^l(n)=s_p^{l+1}(n)\}$.
\end{lem}

\noindent{\bf Proof.}  By definition, $T_p(n)=s_p^j(n)$ where
$j=\min\{ l\geq 0\mbox{ } |\mbox{ } s_p^l(n)=s_p^{l+1}(n)\}$. If
$J\in \{l\geq 0\mbox{ } |\mbox{ } s_p^l(n)=s_p^{l+1}(n)\}$, then
$J\geq j$ and $s_p^j(n)=s_p^J(n)$ thanks to Lemma \ref{Existencia de
Tpn}.
\begin{flushright}$\square$\end{flushright}

\begin{lem}\label{Potencias n y tpn}
For all $x\in \mathbb F_p$ and $n\geq 1$, we have that $x^n=x^{T_p(n)}$.
\end{lem}

\noindent{\bf Proof.} If $n\leq p-1$, then $T_p(n)=s_p^0(n)=n$, so
$x^n=x^{T_p(n)}$. Suppose that $x^m=x^{T_p(m)}$ for all $m<n$ where
$n\geq p$ and we express $n=\sum_{i=0}^s e_ip^i$ in base $p$
expansion where $e_s\neq 0$. Then,
$$
x^n=x^{\sum e_ip^i}=\prod_i x^{e_ip^i}=\prod_i x^{e_i}=x^{\sum
e_i}=x^{s_p(n)}\mod p
$$
where $s_p(n)<n$ because $n\geq p$. By hypothesis, if
$T_p(s_p(n))=s_p^j(s_p(n))$, then
$$
x^n=x^{s_p(n)}=x^{T_p(s_p(n))}=x^{s_p^j(s_p(n))}= x^{s_p^{j+1}(n)}
$$
Observe that
$s_p^{j+1}(n)=s_p^j(s_p(n))=s_p^{j+1}(s_p(n))=s_p^{j+2}(n)$. So,
$T_p(n)=s_p^{j+1}(n)$ by Lemma \ref{Igualdad de Tpn}.
\begin{flushright}$\square$\end{flushright}

\begin{lem}\label{Ecuaciones modulo p}
Let $p$ be a prime. Then, for all $m$ such that $1<m<p$, there
exists a finite number of elements $a_i\in \mathbb F_p^\ast$
(multiplicative group) such that
$$
\left\{\begin{array}{l}
\sum_{i} a_i=1\mod p\\
\sum_{i} a_i^m=0\mod p
\end{array}
\right.
$$
\end{lem}

\noindent{\bf Proof.} Note that $p>2$ because there is not $m\in
\mathbb N$ such that $1<m<2$. Since $\mathbb F_p^\ast$ is a cyclic
group, there exists $g\in \mathbb F^\ast_p$ a generator of $\mathbb
F_p^\ast=\{g,g^2,\ldots, g^{p-1}=1\}$, so $g\neq g^m$ for all
$m=2,\ldots,p-1$. We call $a'_0=g$ and let us consider $h= g^m \mod
p$ with $0<h<p$. Then, we put $a'_i=1$ for $i=1,\ldots, p-h$. In
this case, $$ \sum_{i=0}^{p-h} (a'_i)^m=g^m+\sum_{i=1}^{p-h}
1=g^m+p-h=0\mod p
$$
and
$$
\sum_{i=0}^{p-h} a'_i= g+\sum_{i=1}^{p-h}=g+p-h=g-h\neq 0\mod p
$$
because $h=g^m \mod p$ and, if $g=h\mod p$ then $g=g^m\mod p$!!!. If
we define $a_i=a'_i/(g-h)$, we have the result.
\begin{flushright}$\square$\end{flushright}

\begin{teo}\label{Resultado parcial para los saltos en la
integrabilidad (Tpn)} Let $k$ be a ring of characteristic $p>0$ and
$A$ a $k$-algebra. Let $n\geq 1$ be an integer such that $T_p(n)\neq
1$. Then, $\IDer_k(A;n-1)=\IDer_k(A,n)$.
\end{teo}

\noindent{\bf Proof.} Since $A$ is a $k$-algebra, we can see $A=R/I$
where $R$ is a polynomial ring (in an arbitrary number of variables)
and $I\subseteq R$ an ideal. By Corollary \ref{Saltos del cociente},
$A$ has not leap at $n$ if and only if $\IDer_k(\log
I;n-1)=\IDer_k(\log I;n)$. The inclusion $\IDer_k(\log
I;n-1)\supseteq \IDer_k(\log I;n)$ is always true. Let us consider
$\delta\in \IDer_k(\log I;n-1)$ and  $D\in \HS_k(\log I;n-1)$ an
integral of $\delta$. By Corollary \ref{HSenpolise extiende}, we can
integrate $D$ until $n$. So, we rewrite
$D=(\Id,D_1,\ldots,D_{n-1},D_n)\in \HS_k(R;n)$ as an integral of the
previous $D$ and we obtain an integral of $D_1=\delta$ which is
$(n-1)-I$-logarithmic.

Let us consider $(a_i)_{i}$ a solution of the system of Lemma
\ref{Ecuaciones modulo p} where $m=T_p(n)$. Then,
$$
E:=\circ_i \left(a_i\bullet D\right)= \left(\Id, \sum_i a_i
D_1,\ldots, \sum_i a_i^nD_{n}+ H\right) \mbox{ where }
H:=\sum_{|\beta|=n: \beta_i<n \forall i} \circ_{i}
\left(a_i^{\beta_i}D_{\beta_i}\right)
$$
By Lemma \ref{Potencias n y tpn}, $\sum a_i^n=\sum a_i^{T_p(n)}=0 \mod p$. Moreover, since
$D_\beta(I)\subseteq I$ for all $\beta<n$, $H$ is an $I$-logarithmic differential operator. So, since $\sum a_i=1\mod p$,
$$
E=\left(\Id, D_1,\ldots, H\right)\in\HS_k(\log I;n).
$$
Therefore $D_1\in \IDer_k(\log I;n-1)=\IDer_k(\log I,n)$ and, by
 Corollary \ref{Saltos del cociente}, $\IDer_k(A;n-1)=\IDer_k(A;n)$.
\begin{flushright}$\square$\end{flushright}

\begin{cor}\label{Caso 2p}
Let $k$ be a ring of characteristic $p\geq 3$ and $A$ a
 $k$-algebra. Then, $\IDer_k(A;2p-1)=\IDer_k(A,2p)$.
\end{cor}

\noindent{\bf Proof.} Since $T_p(2p)=2$, we have the result by
Theorem \ref{Resultado parcial para los saltos en la integrabilidad
(Tpn)}. \vspace{-0.68cm}
\begin{flushright}$\square$\end{flushright}

\section{Integrating the first non-vanishing component of a Hasse-Schmidt derivation}

In this section, $k$ will be a commutative ring, $A$ a commutative
$k$-algebra and $I\subseteq A$ an ideal. We start with some
numerical properties that we will use in later sections and we end
up calculating an integral for the first non-vanishing component of
a Hasse-Schmidt derivation which will be the key to prove the main
theorem of section \ref{Ultima seccion}.

\subsection{Numerical results}

In this section, we give some numerical results that will be useful in later results.

\begin{Def}\label{Cpmes}
Let $p,s,m,e$ be integers such that $p,s\geq 1$. Then, we define
$$
C^p_{m,e,s}:=\{j\in \mathbb N \mbox{ }|\mbox{ } mp^j<ep^s\}.
$$
\end{Def}

\begin{lem}\label{Maximo si m entre e y eps}
If $e\leq m<ep^s$, then $C^p_{m,e,s}$ is not empty and  $0\leq \max
C^p_{m,e,s}<s$.
\end{lem}

\noindent{\bf Proof.} $C^p_{m,e,s}\neq\emptyset$ because $j=0$ holds
the inequality, so $\max C^p_{e,m,s}\geq 0$. On the other hand, let
us consider $r\geq s$, then
$$
ep^s\leq ep^r\leq mp^r
$$
so, $r\not\in C^p_{m,e,s}$ and, we can conclude that $0\leq \max
C^p_{m,e,s}<s$.
\begin{flushright}$\square$\end{flushright}

\begin{lem}\label{Desigualdad 1}
Let us assume that $e<m<ep^s$ with $m\neq 0\mod e$ and we denote
$r=\max C^p_{m,e,s}$. Then, $mp^{r+1}-1\geq ep^s$.
\end{lem}

\noindent{\bf Proof.} Since $r=\max C^p_{m,e,s}$, we have that
$mp^{r+1}\geq ep^s$. We will see that the equality never holds.
Suppose that $mp^{r+1}=ep^s$. From Lemma \ref{Maximo si m entre e y
eps}, $r+1\leq s$, so $m=ep^{s-(r+1)}$ but $m$ is not a multiple of
$e$ by hypothesis. Therefore, $mp^{r+1}>ep^s$ and we have the
result.
\begin{flushright}$\square$\end{flushright}

Let us consider $p$ a prime and $n=e_sp^s+\cdots+e_tp^t$ a positive
integer expressed in base $p$ expansion where $s\geq t\geq 1$ and
$0\leq e_i<p$ with $e_s,e_t\neq 0$ (note that $s$ and $t$ can be the
same).

\begin{lem}\label{Max C para N+1}
Let $p$ and $n$ be as before. For all $a\in \mathbb N$ such that
$2p^t+1\leq a<n+1$, we have $0\leq \max C^p_{a,n+1,t}\leq s$.
\end{lem}

\noindent{\bf Proof.} Observe that $0\in C^p_{a,n+1,t}$, so these
sets are not empty. Consider $r>s$, then
$$
(2p^t+1)p^r <(n+1)p^t \Leftrightarrow
(2p^t+1)p^{r-t}=2p^r+p^{r-t}<n+1
$$
The last inequality is false because $n<p^{s+1}\leq p^r$, so
$n+1<p^r+1\leq 2p^r+p^{r-t}$. Hence, $C^p_{2p^t+1,n+1,t}\leq s$.
Now, we consider $a>2p^t+1$ and, as before, $r>s$, then,
$$
(n+1)p^t\leq(2p^t+1)p^r<ap^r
$$
where the first inequality holds  because $\max
C^p_{2p^t+1,n+1,t}\leq s$.  So, $r\not\in C^p_{a,n+1,t}$ for $r>s$,
i.e, $\max C^p_{a,n+1,t}\leq s$.
\begin{flushright}$\square$\end{flushright}

\begin{lem}\label{2p^t<=n}
Let $p,n$ be as before and let us suppose that $n$ is not a power of
$p$, then $2p^t\leq n$.
\end{lem}

\noindent{\bf Proof.} If $p=2$, then $s>t$, because $n$ is not a
power of 2,  so $2p^t=p^{t+1}\leq p^s\leq p^s+e_{s-1}p^{s-1}+\cdots
+p^t=n$. Let us assume that $p\neq 2$. If $s>t$ then, $2p^t<p^s\leq
e_sp^s\leq n$ and we have the inequality. Otherwise, if $s=t$, then
$e_t\geq 2$ because $n$ is not a power of $p$, so $2p^t\leq
e_tp^t=n$.
\begin{flushright}$\square$\end{flushright}

\begin{lem}\label{Desigualdad 2}
Let $p,n$ be as before. Let us assume that $n$ is not a power of $p$
(remember that $n$ is a multiple of $p$). For each integer $a$ such
that $2p^t+1\leq a<n+1$, we denote $r_a=\max C^p_{a,n+1,t}$. Then,
$ap^{r_a+1}-1\geq (n+1)p^t$.
\end{lem}

\noindent{\bf Proof.} By definition $ap^{r_a+1}\geq (n+1)p^t$. We
will see that the equality never holds. Let us suppose that
$ap^{r_a+1}=(n+1)p^t$. Since $a<n+1$, we have that $r_a+1>t$. Then,
$ap^{r_a+1-t}=n+1$, i.e, $n+1$ has to be a multiple of $p$!!! So,
$ap^{r_a+1}-1\geq (n+1)p^t$.
\begin{flushright}$\square$\end{flushright}

\subsection{Computing the integral of the first non-vanishing component of a Hasse-Schmidt derivation}

Let $k$ be a commutative ring and $A$ a commutative $k$-algebra.

\begin{hip}\label{HI}
Let $a\geq 1$ and $p\geq 2$ be integers and $I\subseteq A$ an ideal.
We say that $A$ satisfies the condition $H^I_{p,a}$ if for all $M\in
\mathbb N_+$ not a power of $p$ with $1< M<p^a$, then $\IDer_k(\log
I;M-1)=\IDer_k(\log I;M)$.
\end{hip}

\begin{nota}
Note that $I$ can be $A$. In this case, the condition in Hypothesis \ref{HI} is $\IDer_k(A;M-1)=\IDer_k(A;M)$.
\end{nota}

\begin{lem}\label{La hipotesis HI se cumple}
\begin{enumerate}
\item If $A$ satisfies $H^I_{p,a}$ for some $a\geq 1$, then $A$ satisfies $H^I_{p,s}$ for all $1\leq s\leq a$.
\item If $\charr(k)=p>0$, then $A$ satisfies $H^A_{p,1}$.
\item If $\charr(k)=p>0$ and $A=k[x_i|\mbox{ }i\in \mathcal I]$, the polynomial ring in an arbitrary number of variables, then $A$ satisfies $H^I_{p,1}$ for all $I\subseteq A$ ideal.
\end{enumerate}
\end{lem}

\noindent{\bf Proof.}

\begin{enumerate}
\item It is obvious.
\item If $1<M<p$, then  $M$ can not be a multiple of $p$, i.e., $M\neq 0\mod p$. By Corollary \ref{No tiene saltos en los no multiplos de p}, $\IDer_k(A;M-1)=\IDer_k(A;M)$. 

\item From Corollary \ref{Saltos del cociente},
we have that $\IDer_k(\log I;M-1)=\IDer_k(\log I;M)$ if and only if
$A/I$ does not have leap at $M$, i.e, if
$\IDer_k(A/I;M-1)=\IDer_k(A/I;M)$. Since $A/I$ satisfies
$H^{A/I}_{p,1}$, for all $M\in \mathbb N_+$ with $1<M<p$, we have
the last equality, so $A$ satisfies $H^I_{p,1}$.
\end{enumerate}
\vspace{-0.8cm}\begin{flushright}$\square$\end{flushright}

From now on, $k$ will be a commutative ring, $A$ a commutative $k$-algebra,
$I\subseteq A$ an ideal and $p\geq 2$ an integer.

\begin{lem}\label{a=1, hacer ceros en el siguente tramo}
Let us assume that $A$ satisfies $H^I_{p,1}$. Let $e>1$ be an integer
and $D\in \HS_k(A;ep)$ a Hasse-Schmidt derivation that is
$(ep-1)-I$-logarithmic and $\ell(D;e)=i$ with $0<i<p$. Then, there
exists $D'\in \HS_k(A;ep)$ such that $D'$ is $(ep-1)-I$-logarithmic,
$\ell(D';e)\geq i+1$, $D_\alpha'=D_\alpha$ for all $\alpha\leq ie$
and $D_{ep}'=D_{ep}+H$ where $H$ is an $I$-logarithmic differential
operator.
\end{lem}

\noindent{\bf Proof.} Since $\ell(D;e)=i\geq 1$, from Lemma
\ref{l(D;e) y derivaciones}, we have that $D_{ie+\alpha}\in
\Der_k(\log I)$ for all $\alpha=1,\ldots, e-1$ and, thanks to the
condition $H^I_{p,1}$, we know that all derivations are
$I$-logarithmically $(p-1)$-integrable. Let $D^\alpha\in \HS_k(\log
I;p-1)$ be an integral of $D_{ie+\alpha}$ and consider
$E^{D^\alpha,ie+\alpha}\in \HS_k(\log I;(ie+\alpha)p-1)$, defined in
\ref{EDb}, for all $\alpha=1,\ldots,e-1$. Note that
$(ie+\alpha)p-1>iep\geq ep$, so we can truncate all
these derivations. We denote
$E^\alpha:=\tau_{(ie+\alpha)p-1,ep}\left(E^{D^\alpha,ie+\alpha}\right)$.
Remember that $E^\alpha_{ie+\alpha}=-D_{ie+\alpha}$ and
$\ell(E^\alpha;ie+\alpha)=\lceil ep/ie+\alpha\rceil\geq 2$ for all
$\alpha$  because  $ie+\alpha<(i+1)e\leq ep$.

We denote $E:=E^1\circ E^2\circ \cdots \circ E^{e-1}\in \HS_k(\log I; ep)$.
From Lemma \ref{Composicion de derivaciones y l(D;e)}, $\ell(E)>ie$
and $E_{ie+\alpha}=-D_{ie+\alpha}$ for all $\alpha=1,\ldots, e-1$.
From Lemma \ref{Composicion 1}, $D':=D\circ E\in \HS_k(A;ep)$
satisfies the result because $D_\alpha'=D_\alpha$ for all
$\alpha\leq ie$ and $D'_{ie+\alpha}=D_{ie+\alpha}+E_{ie+\alpha}=0$
for $\alpha=1,\ldots,e-1$, so $\ell(D';e)\geq i+1$. Moreover, from
Lemma \ref{Logaritmico y ultima componente}, $D'$ is
$(ep-1)-I$-logarithmic and $D'_{ep}=D_{ep}+H$ where $H$ is an
$I$-logarithmic differential operator.
\vspace{-0.3cm}\begin{flushright}$\square$\end{flushright}

\begin{lem}\label{a=1, l(D) >=e, De es p-integrable}
Let us assume that $A$ satisfies $H^I_{p,1}$. Let $e\geq 1$ be an
integer and  $D\in \HS_k(A;ep)$ an $(ep-1)-I$-logarithmic
Hasse-Schmidt derivation such that  $\ell(D)\geq e$. Then, $D_e$ is
$p$-integrable and there exists $D'\in \HS_k(A;p)$ an integral of
$D_e$ such that $D'$ is $(p-1)-I$-logarithmic  and $D_p'=D_{ep}+H$
where $H$ is an $I$-logarithmic differential operator.
\end{lem}

\noindent{\bf Proof.} If $e=1$, the result is trivial, so we will
suppose that $e>1$. We prove the result by reverse induction on
$p\geq \ell(D;e)\geq 1$ (note that $\ell(D;e)\geq 1$ because
$\ell(D)\geq e$ and $\ell(D;e)\leq \lceil ep/e\rceil=p$).

If $\ell(D;e)=p$, by Lemma \ref{Operacion inversa de desplazar},
there exists $D'\in \HS_k(A;p)$ such that $D_\alpha'=D_{\alpha e}$ for all $\alpha=1,\ldots,p$.
Since $D$ is $(ep-1)-I$-logarithmic, $D'_\alpha(I)= D_{\alpha
e}(I)\subseteq I$ for all $\alpha<p$, so $D'$ is
$(p-1)-I$-logarithmic. Moreover, $D'_p=D_{ep}$ so, $D'$ satisfies the lemma.

Now, let us assume that if $D$ is a Hasse-Schmidt derivation with
$\ell(D)\geq e$ and $\ell(D;e)\geq i+1$ where $1\leq i<p$, then we
have the result. We will prove it for a Hasse-Schmidt derivation $D$
such that $\ell(D)\geq e$ and $\ell(D;e)=i$.

By Lemma \ref{a=1, hacer ceros en el siguente tramo}, there exists
$D'\in \HS_k(A;ep)$ such that $D'$ is $(ep-1)-I$-logarithmic,
$\ell(D';e)\geq i+1$, $D'_\alpha=D_\alpha$ for all $\alpha\leq ie$
and $D_{ep}'=D_{ep}+H$ where $H$ is an $I$-logarithmic differential operator. Since $\ell(D')\geq
e$ because $\ell(D';e)\geq i+1\geq 1$, we can apply the induction hypothesis, so
there exists $D''\in \HS_k(A;p)$ an integral of $D_e'=D_e$ which is $(p-1)-I$-logarithmic and
$D_p''=D_{ep}'+H'=D_{ep}+H+H'$ where $H+H'$ is an $I$-logarithmic differential operator. Hence, we
have the lemma.
\vspace{-0.3cm}\begin{flushright}$\square$\end{flushright}

\begin{lem}\label{Extension de integral de -D1}
Let us assume that $A$ satisfies $H^I_{p,a}$ for some $a\geq 1$. Let
$e,s,m$ be integers such that $1\leq s\leq a$ and $1< e\leq m<
ep^s$. We denote $r:=\max C^p_{m,e,s}$ and we consider $\delta\in
\IDer_k(\log I;p^r)$. We have the following properties:
\begin{itemize}
\item[1] If $m=0\mod e$, then there exists $E\in \HS_k(\log I;ep^s-1)$ such that $E_m=-\delta$ and $\ell(E;m)=\lceil ep^s-1/m\rceil$.
\item[2] If $m\neq 0\mod e$, then there exists $E\in \HS_k(\log I;ep^s)$ such that $E_m=-\delta$ and $\ell(E;m)=\lceil ep^s/m\rceil$.
\end{itemize}
\end{lem}

\noindent{\bf Proof.} By Lemma \ref{Maximo si m entre e y eps}, we
have that $0\leq r<s\leq a$, so $p^{r+1}\leq p^a$. Thanks to the
condition $H^I_{p,a}$, we have that $\delta\in \IDer_k(\log
I;p^r)=\IDer_k(\log I; p^{r+1}-1)$. Let $D\in \HS_k(\log
I;p^{r+1}-1)$ be an integral of $\delta$. Let us consider
$E^{D,m}\in \HS_k(\log I,mp^{r+1}-1)$ (Definition \ref{EDb}) where
$E^{D,m}_m=-\delta$ and $\ell\left(E^{D,m};m\right)=p^{r+1}$, i.e., $E^{D,m}_\alpha=0$ for all $\alpha\neq
0\mod m$. Moreover, from the definition of $r$, $mp^{r+1}-1\geq
ep^s-1$.

If $m=0\mod e$, then $E=\tau_{mp^{r+1}-1,ep^s-1}(E^{D,m})$ satisfies
the lemma. Otherwise, if $m\neq 0\mod e$, by Lemma \ref{Desigualdad
1}, $mp^{r+1}-1\geq ep^s$.  So, $E=\tau_{mp^{r+1}-1,ep^s}(E^{D,m})$
satisfies the lemma.
\begin{flushright}$\square$\end{flushright}

\begin{lem}\label{Quitar la primera comp no 0 y dejar las m-1 siguientes =}
Let us assume that $A$ satisfies $H^I_{p,a}$ for some $a\geq 1$. Let
$e,s,m$ be integers such that $1\leq s\leq a$ and $1<e\leq m< ep^s$
and we denote $r:=\max C^p_{m,e,s}$. Let $D\in \HS_k(\log I;ep^s-1)$
be a Hasse-Schmidt derivation such that $\ell(D)\geq m$ where
$D_m\in \IDer_k(\log I;p^r)$. Then, there exists $D'\in \HS_k(\log
I;ep^s-1)$ such that $\ell(D')\geq m+1$ with $D'_\alpha=D_\alpha$
for all $\alpha=m+1,\ldots, 2m-1$.
\end{lem}

\noindent{\bf Proof.} If $D_m=0$, we put $D'=D$ and we have the
lemma. Let us assume that $D_m\neq 0$. From Lemma \ref{Extension de
integral de -D1}, for any $e\leq m<ep^s$, we have $E\in \HS_k(\log
I;ep^s-1)$ such that $E_m=-D_m$ and $\ell(E;m)=\lceil
ep^s-1/m\rceil$. We can apply Lemma \ref{composicion 2} to
$D'=D\circ E\in \HS_k(\log I; ep^s-1)$. 
Then, $\ell(D')\geq m$ and
$$
D'_\alpha=\left\{\begin{array}{ll}
D_m+E_m&\mbox{ if } \alpha=m\\
D_\alpha&\mbox{ if }\alpha=m+1,\ldots,2m-1
\end{array}
\right.
$$
Since $E_m=-D_m$, $D'_m=0$ and hence, $\ell(D')\geq m+1$ and $D'$
satisfies the lemma.
\begin{flushright}$\square$\end{flushright}

\begin{teo}\label{Extraer una integral de De a partir de D}
Let us suppose that $A$ satisfies $H^I_{p,a}$ for some $a\geq 1$.
Let $e,s\geq 1$ be two integers such that $s\leq a$ and $D\in
\HS_k(A;ep^s)$ an $(ep^s-1)-I$-logarithmic Hasse-Schmidt derivation
with $\ell(D)\geq e$. Then, there exists $D'\in \HS_k(A;p^s)$ an
integral of $D_e$ such that $D'$ is $(p^s-1)-I$-logarithmic with
$D_{p^s}'=D_{ep^s}+H$ where $H$ is an $I$-logarithmic differential
operator.
\end{teo}

\noindent{\bf Proof.} We prove the result by induction on $s\geq 1$.
Note that if $s=1$, we have the theorem from Lemma \ref{a=1, l(D)
>=e, De es p-integrable}. So, let us assume that the theorem is true
for all $j$ such that $1\leq j< s\leq a$. Note that we can suppose
that $e>1$ (if $e=1$ the theorem is trivial).  We will divide this
proof in several lemmas:

\begin{itemize}
\item[]

\begin{lem}\label{p^r-integrable}
Let $D\in \HS_k(\log I;ep^s-1)$ such that $\ell(D)\geq m$ with
$ep^s>m\geq e>1$. Then, $D_m\in \IDer_k(\log I;p^r)$ with $r=\max
C^p_{m,e,s}<s$.
\end{lem}

\noindent{\bf Proof.} By Lemma \ref{Maximo si m entre e y eps}, we
have that $0\leq r<s$. We rewrite $D:=\tau_{ep^s-1,mp^r}(D)\in
\HS_k(\log I;mp^r)$ (note that $mp^r\leq ep^s-1$ by definition of
$C^p_{m,e,s}$). If $r=0$, then it is obvious that $D_m$ is $I$-logarithmically $p^r$-integrable. Let us suppose that $r\geq 1$.  Then, since $1\leq r<s\leq a$, by the induction hypothesis of
the theorem, there exists $D'\in \HS_k(A;p^r)$ an integral of $D_m$
such that $D'$ is $(p^r-1)-I$-logarithmic and
$D'_{p^r}=D_{mp^r}+(\mbox{\it some $I$-logarithmic diff. op.)}$. But
$D_{mp^r}$ is $I$-logarithmic, so $D'$ is $I$-logarithmic too and
$D_m$ is $I$-logarithmically $p^r$-integrable.
\begin{flushright}$\diamond$\end{flushright}


\begin{lem}\label{Quitar las i primeras comp no 0 y dejar las e-1 siguientes}
Let $D\in \HS_k(\log I;ep^s-1)$ such that $\ell(D;e)=i<p^s$ and
$\ell(D)\geq e>1$. Then, there exists $D'\in \HS_k(\log I;ep^s-1)$
such that $\ell(D')>ie$ and $D'_\alpha=D_\alpha$ for all
$\alpha=ie+1,\ldots,ie+e-1$.
\end{lem}

\noindent{\bf Proof.} Note that the only components that can be
not zero before $ie+1$ are those that are in the multiples of $e$.
If $\ell(D)>ie$ then the lemma is obvious, otherwise $\ell(D)=je$
for some $1\leq j\leq i$. We  will prove the result by reverse
induction on $1\leq j\leq i$.

Let us assume that $\ell(D)= ie$. By Lemma \ref{p^r-integrable},
$D_{ie}\in \IDer_k(\log I;p^r)$ where $r=\max C^p_{ie,e,s}<s$. From
Lemma \ref{Quitar la primera comp no 0 y dejar las m-1 siguientes
=}, there exists $D'\in \HS_k(\log I;ep^s-1)$ such that
$\ell(D')\geq ie+1$ and $D_\alpha'=D_\alpha$ for all
$\alpha=ie+1,\ldots,\min\{ ep^s-1,2ie-1\}$. Note that $ie+e-1\leq
ep^s-1$, so $D'$ satisfies the lemma.

Let us suppose now that the lemma is true for all derivation with
$\ell(\cdot)>je$ and we will prove it for $j<i$.

By Lemma \ref{p^r-integrable}, $D_{je}\in \IDer_k(\log I;p^r)$ where
$r=\max C^p_{je,e,s}<s$. From Lemma \ref{Extension de integral de
-D1}, there exists $E\in \HS_k(\log I;ep^s-1)$ such that
$E_{je}=-D_{je}$ and $\ell(E;je)=\lceil ep^s-1/je\rceil\geq 1$. We
can apply Lemma \ref{composicion 2} to $D$ and $E$ and we obtain
$D'=D\circ E\in \HS_k(\log I,ep^s-1)$ such that $\ell(D')\geq je$,
$\ell(D';e)\geq \ell(D;e)=i$ and
$$ D'_\alpha=\left\{\begin{array}{ll}
D_{je}+E_{je}&\mbox{if }\alpha=je\\
D_\alpha&\mbox{for all } \alpha=ie+1,\ldots,ie+e-1
\end{array}
\right.
$$
Since $\ell(D;e)=i$, there exists $a\in \{1,\ldots,e-1\}$ such that
$D_{ie+a}\neq0$ and, since $D'_{ie+a}=D_{ie+a}$, we have that
$\ell(D';e)=i$. Moreover, $E_{je}=-D_{je}$, so $\ell(D')\geq je+1$,
but $\ell(D';e)>j$, therefore $\ell(D')\geq (j+1)e$. Now, we can
apply the induction hypothesis. Hence, there exists $D''\in
\HS_k(\log I;ep^s-1)$ such that $\ell(D'')>ie$ and
$D''_\alpha=D'_\alpha=D_\alpha$ for all $\alpha=ie+1,\ldots,ie+e-1$
and we have the lemma.
\begin{flushright}$\diamond$\end{flushright}

\begin{lem}\label{Integrales para los tramos no multiplos de e}
Let $D\in \HS_k(\log I;ep^s-1)$ be a Hasse-Schmidt derivation such
that $\ell(D)>ie$ with $1\leq i<p^s$. Then, for all
$\alpha=1,\ldots, e-1$ there exists $E^\alpha\in \HS_k(\log I;ep^s)$
such that $E^\alpha_{ie+\alpha}=-D_{ie+\alpha}$ and
$\ell(E^\alpha;ie+\alpha)=\lceil ep^s/ie+\alpha\rceil$.
\end{lem}

\noindent{\bf Proof.} If $\ell(D)\geq (i+1)e$, then $D_{ie+\alpha}=0$
for all $\alpha=1,\ldots, e-1$ and we have the result. Let us
suppose that $\ell(D)=(i+1)e-1$. By Lemma \ref{p^r-integrable},
$D_{(i+1)e-1}\in \IDer_k(\log I;p^{r_{e-1}})$ where $r_{e-1}=\max
C^p_{(i+1)e-1,e,s}$. Since $(i+1)e-1\neq 0\mod e$, Lemma
\ref{Extension de integral de -D1} give us the result.

Let us assume that the lemma is true for all Hasse-Schmidt
derivations such that $\ell(\cdot)=ie+\beta$ with $1\leq j<\beta
\leq e-1$ and we will prove it for a Hasse-Schmidt derivation $D$
such that $\ell(D)=ie+j$.

As before, from Lemma \ref{p^r-integrable} and Lemma \ref{Extension
de integral de -D1}, there exists $E^j\in \HS_k(\log I;ep^s)$ such
that $E^j_{ie+j}=-D_{ie+j}$ and $\ell(E^j;ie+j)=\lceil ep^s/ie+j\rceil$.
We can apply Lemma \ref{composicion 2} to $D$ and
$E:=\tau_{ep^s,ep^s-1}(E^j)$ obtaining $D'=D\circ E\in \HS_k(\log
I;ep^s-1)$ such that $\ell(D')\geq ie+j$ and
$$
D_\alpha'=\left\{\begin{array}{ll}
D_\alpha+E_\alpha&\mbox{if }\alpha=ie+j\\
D_\alpha&\mbox{if } \alpha=ie+j+1,\ldots,\min\{ep^s-1,2(ie+j)-1\}
\end{array}
\right.
$$
Note that $ie+e-1\leq ep^s-1$, so $D'_\alpha=D_\alpha$ for all
$\alpha=ie+j+1,\ldots, ie+e-1$. Since $E_{ie+j}=-D_{ie+j}$,
$\ell(D')>ie+j$ and we can use the induction hypothesis on $D'$
obtaining that, for all $\alpha=j+1, \ldots, e-1$, there exists
$E^\alpha\in \HS_k(\log I;ep^s)$ such that
$E_{ie+\alpha}=-D'_{ie+\alpha}=-D_{ie+\alpha}$ and
$\ell(E^\alpha)=\lceil ep^s/ie+\alpha\rceil$. So, we have the lemma.
\begin{flushright}$\diamond$\end{flushright}

\begin{lem}\label{Hacer ceros en el siguiente tramo}
Let $D\in \HS_k(A;ep^s)$ be an $(ep^s-1)-I$-logarithmic Hasse-Schmidt
derivation with $1\leq \ell(D;e)=i<p^s$. Then, there exists $D'\in
\HS_k(A;ep^s)$ such that $D'$ is $(ep^s-1)-I$-logarithmic,
$\ell(D';e)\geq i+1$, $D_{je}'=D_{je}$ for all $j\leq i$ and
$D'_{ep^s}=D_{ep^s}+H$ where $H$ is an $I$-logarithmic differential
operator.
\end{lem}

\noindent{\bf Proof.} Let us consider
$D^\tau=\tau_{ep^s,ep^s-1}(D)\in \HS_k(\log I;ep^s-1)$. Then,
$\ell(D^\tau;e)=i\geq 1$ so, $\ell(D^\tau)\geq e$. By Lemma
\ref{Quitar las i primeras comp no 0 y dejar las e-1 siguientes},
there exists $D'\in \HS_k(\log I;ep^s-1)$ such that $\ell(D')>ie$
and $D'_{ie+\alpha}=D_{ie+\alpha}$ for all $\alpha=1,\ldots, e-1$.

By Lemma \ref{Integrales para los tramos no multiplos de e}, for each
$\alpha=1,\ldots,e-1$, there exists $E^\alpha\in \HS_k(\log I;ep^s)$
such that $E^\alpha_{ie+\alpha}=-D'_{ie+\alpha}=-D_{ie+\alpha}$ and
$\ell(E^\alpha;ie+\alpha)=\lceil ep^s/ie+\alpha\rceil$. Note that
$ie+\alpha<ep^s$ so, $\lceil ep^s/ie+\alpha\rceil\geq 2$. By Lemma
\ref{Composicion de derivaciones y l(D;e)}, if we denote
$E=E^1\circ \cdots \circ E^{e-1}\in \HS_k(\log I;ep^s)$, then
$\ell(E)\geq ie+1$ and
$E_{ie+\alpha}=E^\alpha_{ie+\alpha}=-D_{ie+\alpha}$.

Now, we consider $D'=D\circ E\in \HS_k(A;ep^s)$. By Lemma
\ref{Logaritmico y ultima componente}, $D'$ is $(ep^s-1)-I$-logarithmic and
$D'_{ep^s}=D_{ep^s}+H$ where $H$ is an $I$-logarithmic differential operator. On the other hand,
by Lemma \ref{Composicion 1}, we have that
$$
D'_\beta=\left\{\begin{array}{ll}
D_\beta&\mbox{ if } \beta\leq ie\\
D_{\beta}+E_{\beta}& \mbox{ if } \beta=ie+1,\ldots, ie+e-1
\end{array}
\right.
$$
Hence, $D'_\beta=0$ for all $\beta=ie+1,\ldots,ie+e-1$ so,
$\ell(D';e)\geq i+1$. Therefore, $D'$ satisfies the lemma.
\begin{flushright}$\diamond$\end{flushright}
\end{itemize}

Now, with the help of the previous lemmas we will finish the proof
of Theorem \ref{Extraer una integral de De a partir de D}. We show
this result by reverse induction on $1\leq \ell(D;e)\leq p^s$.

If $\ell(D;e)=p^s$, by Lemma \ref{Operacion inversa de desplazar},
there is $D'\in \HS_k(A;p^s)$ such that $D'_\alpha=D_{\alpha e}$ for
all $\alpha\leq p^s$ and we have the result.

Let us assume that the theorem is true for Hasse-Schmidt derivation
with $\ell(\cdot; e)>i$ for $1\leq i<p^s$ and we will prove it for
$D\in \HS_k(A;ep^s)$ with $\ell(D;e)=i$. By Lemma \ref{Hacer ceros
en el siguiente tramo}, there exists $D'\in \HS_k(A;ep^s)$ which is
$(ep^s-1)-I$-logarithmic, such that $\ell(D';e)\geq i+1$ (so $\ell(D')\geq e$),
$D'_{ej}=D_{ej}$ for all $j\leq i$ and $D'_{ep^s}=D_{ep^s}+ H$ where $H$ is an $I$-logarithmic differential operator.
By induction hypothesis, $D_e'=D_e$ has a
$(p^s-1)-I$-logarithmic $p^s$-integral, $D''\in \HS_k(A;p^s)$ with
$D''_{p^s}=D'_{ep^s}+H'=D_{ep^s}+H+H'$ where $H+H'$ is an $I$-logarithmic differential operator. Hence, we have the result.
\begin{flushright}$\square$\end{flushright}

\section{Integrability and leaps}\label{Ultima seccion}

In this section, we prove that, any commutative $k$-algebra, where $k$ is a commutative
ring of characteristic $p>0$, only has leaps at powers of $p$, or
what is the same:

\begin{teo}\label{Teo Integrabilidad}
Let $k$ be a ring of $\charr(k)=p>0$ and $A$ a
$k$-algebra. Then, for all $n>1$ not a power of $p$,
$\IDer_k(A;n-1)=\IDer_k(A;n)$.
\end{teo}

\noindent{\bf Proof.}  It is enough to show the theorem when $n$ is
a multiple of $p$, not a power of $p$ because, if $n\neq 0\mod p$, by Corollary \ref{No tiene saltos en los no multiplos de p},
we have the result. We will prove this theorem by induction on $n$
multiple of $p$, not a  power of $p$. We have two base cases, when
$p=2$ and $p\neq 2$. In the first case, we have to prove that $\IDer_k(A;5)=\IDer_k(A;6)$. Proposition \ref{Caso p=2}
gives us the theorem. In the second one, we have to prove that $\IDer_k(A;2p-1)=\IDer_k(A;2p)$ and we have the result by
Corollary \ref{Caso 2p}. Let us assume that for all $m<n$ not a
power of $p$, $\IDer_k(A;m-1)=\IDer_k(A;m)$ and we will prove the
equality for $n$, a multiple of $p$, not a power of $p$.

Since $A$ is a $k$-algebra, we can express $A=R/I$ where
$R=k[x_i|\mbox{ }i\in \I]$ is a polynomial ring of an arbitrary number of
variables and $I\subseteq R$ an ideal. Then, by Corollary
\ref{Saltos del cociente}, we have that $\IDer_k(\log
I;m-1)=\IDer_k(\log I;m)$ for all $m<n$ not a power of $p$ and it is
enough to prove that $\IDer_k(\log I;n-1)=\IDer_k(\log I;n)$.

Let us express $n=e_sp^s+\cdots+e_tp^t$ in base $p$ expansion where
$1\leq t\leq s$ and  $0\leq e_i<p$ with $e_s,e_t\neq 0$. By induction hypothesis, we have that $R$ satisfies $H_{p,s}^I$
(\ref{HI}).

Let $\delta\in \IDer_k(\log I;n-1)$ be a $k$-derivation and $D\in
\HS_k(\log I;n-1)$ an integral of $\delta$. We can integrate $D$
until length infinity (Corollary \ref{HSenpolise extiende}), so we rewrite
$D\in \HS_k(R)$ the integral of $D$. Note that $D_1=\delta$ and $D$
is $(n-1)-I$-logarithmic. Now, we consider $G:=G^{D,p^t}\in
\HS_k(R;(n+1)p^t)$ the Hasse-Schmidt derivation defined in
\ref{Definicion de Gpt}. From Lemma \ref{I-logaritmico para
derivacion Gpt}, $G$ is $((n+1)p^t-1)-I$-logarithmic, $\ell(G)\geq
2p^t+1$ and $G_{(n+1)p^t}=\binom{n}{p^t} D_n+H$ where $H$ is an $I$-logarithmic differential operator.

By Lemma \ref{2p^t<=n}, we have that $2p^t+1\leq n+1$. If
$n+1=2p^t+1$, then from Theorem \ref{Extraer una integral de De a
partir de D}, we obtain $T\in
\HS_k(R;p^t)$ a $(p^t-1)-I$-logarithmic Hasse-Schmidt derivation
such that $T_{p^t}=G_{(n+1)p^t}+H'=\binom{n}{p^t} D_n+H+H'$ where $H+H'$ is an $I$-logarithmic differential operator.

Let us suppose now that $2p^t+1<n+1$ and we denote $r=\max
C^p_{2p^t+1,n+1,t}$. By Lemma \ref{Max C para N+1}, $0\leq r\leq s$  and by definition of $C^p_{2p^t+1,n+1,t}$,  $(2p^t+1)p^r<(n+1)p^t$.
Hence, we can consider $\tau_{(n+1)p^t,(2p^t+1)p^r}(G)\in \HS_k(\log
I;(2p^t+1)p^r)$.  If $r=0$, then $G_{2p^t+1}\in \IDer_k(\log I;p^r)$. Otherwise, $r\geq 1$ and applying Theorem \ref{Extraer una integral de De a
partir de D} to this Hasse-Schmidt derivation, we have that
$G_{2p^t+1}$ is $p^r$-integrable, and there exists $D'\in
\HS_k(R;p^r)$ an integral that is $(p^r-1)-I$-logarithmic with
$D'_{p^r}=G_{(2p^t+1)p^r}+\mbox{\it  some $I$ -logarithmic diff. op.}$.
Since $G_{(2p^t+1)p^r}$ is $I$-logarithmic, $G_{2p^t+1}$ is $I$-logarithmically $p^r$-integrable. So, in both cases, we have that $G_{2p^t+1}\in
\IDer_k(\log I;p^r)$.  We have two cases:
\begin{itemize}
\item If $r<s$, then $G_{2p^t+1}\in \IDer_k(\log I;p^{r+1}-1)$ from the hypothesis, i.e, there exists  $D'\in \HS_k(\log I; p^{r+1}-1)$ an integral of
$G_{2p^t+1}$ and we can consider $E^{D',2p^t+1}\in \HS_k\left(\log
I;(2p^t+1)p^{r+1}-1\right)$.

By Lemma \ref{Desigualdad 2}, we can consider $T=
\tau_{(2p^t+1)p^{r+1}-1,(n+1)p^t}\left(E^{D',2p^t+1}\right) \in
\HS_k(\log I;(n+1)p^t)$ where $T_{2p^t+1}=-G_{2p^t+1}$ and
$\ell(T)\geq 2p^t+1$.

\item If $r=s$, then $G_{2p^t+1}\in \IDer_k(\log I;p^s)$.
Since $p^s<n<p^{s+1}$, $G_{2p^t+1}\in \IDer_k(\log I;n-1)$. Let
$D'\in \HS_k(\log I;n-1)$ be an integral of $G_{2p^t+1}$ and let us
consider $E^{D';2p^t+1}\in \HS_k(\log  I;(2p^t+1)n-1)$. Note that
    $$
    n(2p^t+1)-1>(n+1)p^t \Leftrightarrow np^t+n-1>p^t
    $$
Since the last inequality always holds, we can consider
$T=\tau_{(2p^t+1)n-1,(n+1)p^t)}\left(E^{D',2p^t+1}\right)\in
\HS_k(\log I;(n+1)p^t)$ where $T_{2p^t+1}=-G_{2p^t+1}$ and
$\ell(T)\geq 2p^t+1$.
\end{itemize}

Therefore, in both cases, we can compose $G$ and $T$ obtaining
$$
G^{(1)}:=T\circ
G^{p^t}=\left(\Id,0,\ldots,0,G_{2p^t+2}^{(1)},\ldots,
G_{n+1}^{(1)},\ldots, G_{(n+1)p^t}^{(1)}\right)\in \HS_k(R;(n+1)p^t)
$$
a $((n+1)p^t-1)-I$-logarithmic Hasse-Schmidt derivation where
$G^{(1)}_{(n+1)p^t}=\binom{n}{p^t} D_n+H$ with $H$ an
$I$-logarithmic differential operator.

Suppose that, by doing this process, we obtain a
$((n+1)p^t-1)-I$-logarithmic Hasse-Schmidt derivation:
$$
G^{(j)}=(\Id,0,\ldots,0,G_j^{(j)},\ldots,
G_{n+1}^{(j)},\ldots,\binom{n}{p^t} D_n+H)\in
\HS_k\left(R;(n+1)p^t\right)
$$
with $H(I)\subseteq I$ and $2p^t+1<j<n+1$.

We denote $r=\max C^p_{j,n+1,t}$. By Lemma \ref{Max C para N+1},
$0\leq r\leq s$. Since $jp^r<(n+1)p^t$, we have
$\tau_{(n+1)p^t,jp^r}\left(G^{(j)}\right)\in \HS_k(\log I;jp^r)$ and
we can deduce that $G_j^{(j)}$ is $I$-logarithmically
$p^r$-integrable in the same way as before. We have two cases:
\begin{itemize}
\item If $r<s$, then $G^{(j)}_{j}\in \IDer_k(\log I;p^r)=\IDer_k(\log I;p^{r+1}-1)$. Let us consider $D'\in \HS_k(\log I;p^{r+1}-1)$ an integral of $G^{(j)}_j$ and $E^{D',j}\in \HS_k(\log I;jp^{r+1}-1)$. By Lemma \ref{Desigualdad 2}, $jp^{r+1}-1\leq (n+1)p^t$. So, we have $T=\tau_{jp^{r+1}-1,(n+1)p^t}(E^{D',j})\in \HS_k(\log I;(n+1)p^t)$ where $T_j=-G^{(j)}_j$ and $\ell(T)\geq j$.

\item If $r=s$, then $G^{(j)}_j\in \IDer_k(\log I;p^s)=\IDer_k(\log I;n-1)$. Then, there exists $D'\in \HS_k(\log I;n-1)$ an integral of $G^{(j)}_j$ and we can consider $E^{D',j}\in \HS_k(\log I;jn-1)$. Since $jn-1>(2p^t+1)n-1>(n+1)p^t$, we can define $T=\tau_{jn-1,(n+1)p^t}(E^{D',j})\in \HS_k(\log I;(n+1)p^t)$ where $T_j=-G^{(j)}_j$ and $\ell(T)\geq j$.
\end{itemize}
Therefore, we can obtain a $((n+1)p^t-1)-I$-logarithmic Hasse-Schmidt derivation:
$$
G^{(j+1)}:=T\circ G^{(j)}=G^{(j+1)}=
(\Id,0,\ldots,0,G_{j+1}^{(j+1)},\ldots,G^{(j+1)}_{n+1},\ldots,
\binom{n}{p^t}D_n+H')\in \HS_k\left(R;(n+1)p^t\right)
$$
where $H'$ is an $I$-logarithmic differential operator.
So, we can do this process for all $j$ such that $2p^t+1\leq j<n+1$
and we obtain a $((n+1)p^t-1)-I$-logarithmic Hasse-Schmidt derivation:
$$
G^{(n+1)}=(\Id,0,\ldots,0, G^{(n+1)}_{n+1},\ldots,
\binom{n}{p^t}D_n+H')\in \HS_k\left(R;(n+1)p^t\right)
$$
where $H'$ is an $I$-logarithmic differential operator. Then, we
can apply Theorem \ref{Extraer una integral de De a partir de D} to
$G^{(n+1)}$. So, in both cases, when $n+1=2p^t+1$ or not, we have
that there exists $T=(\Id,T_1,\ldots,\binom{n}{p^t}D_n+H')\in
\HS_k(R;p^t)$ that is $(p^t-1)-I$-logarithmic and $H'$ is an $I$-logarithmic differential operator.

Let $f\in \mathbb F_p^\ast$ be the inverse of $\binom{n}{p^t}$. So
that,
$$
\begin{array}{cc}
\displaystyle D\circ \left(-f\bullet T\right)[n/p^t]=\left(\Id,D_1,
\ldots,
D_n+(-f)^{p^t}\binom{n}{p^t}D_n-f^{p^t}H'+
\sum_{\substack{\alpha+\beta=n \\ \alpha,\beta\neq
0}}D_\alpha\circ \left(\left(-f\cdot T\right)[n/p^t]\right)_\beta
\right)=\\
\displaystyle
\left(\Id,D_1,\ldots,\sum_{\substack{\alpha+\beta=n \\ \alpha,\beta\neq
0}}
D_\alpha\circ \left(\left(-f\cdot
T\right)[n/p^t]\right)_\beta-fH'\right)\in \HS_k(\log I;n)
\end{array}
$$
Hence, $D_1=\delta\in\IDer_k(\log I;n)$.
\begin{flushright}$\square$\end{flushright}

{\bf Acknowledgment.} The author thanks Professor Luis Narv\'aez
Macarro for his careful reading of this paper with numerous useful
comments.

\end{document}